\documentclass[twoside,12pt]{article}
\usepackage[]{amsmath,amsthm,amssymb}
\usepackage[latin1]{inputenc}

\begin{document}

\title{{\Large\bf  New biorthogonal sequences generated by index integrals of the weight functions}}

\author{Semyon   YAKUBOVICH}
\maketitle

\markboth{\rm \centerline{ Semyon YAKUBOVICH}}{}
\markright{\rm \centerline{Orthogonal polynomials via index integrals}}

\begin{abstract} {We exhibit new biorthogonal sequences generated by index integrals of the squares of the modified Bessel functions and gamma functions. The composition orthogonality, involving differential operators is employed. Generalized Wilson polynomials are introduced. Some properties are investigated.}

\end{abstract}
\vspace{4mm}



\vspace{4mm}

\section {Introduction and preliminary results}

The aim of the present contribution is to investigate properties of biorthogonal sequences with respect to $L_2$-inner product in terms of the index integrals, involving squares of the modified Bessel functions and Euler's  gamma functions  [3].  We emphasize the role of the composition orthogonality and hypergeometric orthogonal polynomials  [2].  Recent attempts of such investigation is done by the author in [9], [10].  In fact, considering the kernel of the Kontorovich-Lebedev transform [6] $K_{i\tau}(x)$ as the modified Bessel or Macdonald function of the pure imaginary index $i\tau,\ \tau \in \mathbb{R}$ and positive argument $x >0$. This is an even with respect to $\tau$ and  real-valued  function, and we are motivated by its index Fourier cosine transform  (see Entry 2.16.48.19 in [4, Vol. II])

$$\int_0^\infty \cos(\tau y)K_{i\tau}(x) d\tau = {\pi\over 2}\ e^{-x\cosh y},\quad x >0,\ y \in \mathbb{R}.\eqno(1.1)$$
An immediate consequence of this formula is the termwise differentiation by $y$ under the integral sign, which can be easily justified by virtue of the absolute and uniform convergence. As a result we arrive at the following relation

$$  \int_0^\infty \cos(\tau y) \tau^{2n} K_{i\tau}(x) d\tau =  (-1)^n {\pi\over 2}\  {\partial^{2n}\over \partial y^{2n}} \bigg[ e^{-x\cosh y}\bigg],\quad n \in \mathbb{N}_0.\eqno(1.2)$$ 
Moreover, taking the Mellin transform by $x$ [5] through  (1.1), we interchange the order of integration, owing to the dominated convergence which is based on the uniform inequality for the Macdonald function (see [6, formula (1.100)])

$$\left|K_{\nu+i\tau}(x)\right| \le e^{-\delta |\tau|} K_\nu(\cos\delta),\quad x >0,\ \tau, \nu \in \mathbb{R},\ \delta \in \left[0,\ {\pi\over 2}\right).\eqno(1.3)$$
Then,  employing  relation 2.16.2.2 in [4, Vol. II] to derive formula (1.104) in [6] 

$$\int_0^\infty \cos(\tau y) \Gamma\left({s+i\tau\over 2}\right) \Gamma\left({s-i\tau\over 2}\right) d\tau = {\pi\over 2^{s-1}}\ {\Gamma(s)\over \cosh^s y},\quad {\rm Re} s >0,\ y \in \mathbb{R},$$
where $\Gamma (z)$ is Euler's gamma function [3].  Again, differentiating termwisely  with respect to $y$ through the latter equality under the same justification, we deduce the representation

$$\int_0^\infty \cos(\tau y) \tau^{2n} \Gamma\left({s+i\tau\over 2}\right) \Gamma\left({s-i\tau\over 2}\right) d\tau$$

$$ = (-1)^n  {\pi\over 2^{s-1}} \Gamma(s)  {\partial^{2n}\over \partial y^{2n}} \bigg[ {1\over \cosh^s y}\bigg],\ n \in \mathbb{N}_0.\eqno(1.4)$$
It is known [3]  that the modified Bessel function $K_{i\tau}(x)$ is an eigenfunction of the following second order differential operator
$${\cal A} \equiv x^2- x{d\over dx} x {d\over dx},\eqno(1.5)$$
i.e. we have
$${\cal A} \ K_{i\tau}(x)= \tau^2 K_{i\tau}(x).\eqno(1.6)$$
Hence, using (1.1) and the differentiation under the integral sign, we get from (1.2) the following operator equality

$${\cal A}^n \bigg[ e^{-x\cosh y} \bigg] = (-1)^n  D_y^{2n} \bigg[ e^{-x\cosh y}\bigg],\quad n \in \mathbb{N}_0,\eqno(1.7)$$  
where ${\cal A}^n$ is the $n$th iteration of the operator (1.5) and $D_y\equiv {\partial\over \partial y}$.  It can be extended to any polynomial $q_n(x)$ of degree at most $n$ to obtain

$$q_n\left( {\cal A}\right) \bigg[ e^{-x\cosh y} \bigg] =  q_n\left( - D_y^2\right) \bigg[ e^{-x\cosh y} \bigg].\eqno(1.8)$$
Meanwhile, extending equality (1.4) to an arbitrary polynomial $q_n$, it reads

$$\int_0^\infty \cos(\tau y) q_n(\tau^{2})  \Gamma\left({s+i\tau\over 2}\right) \Gamma\left({s-i\tau\over 2}\right) d\tau$$

$$ =   {\pi\ \Gamma(s) \over 2^{s-1}} \  q_n\left( - D_y^2\right) \bigg[  {1\over \cosh^s y} \bigg],\quad n \in \mathbb{N}_0.\eqno(1.9)$$ 
Further, we employ the product formula for the modified Bessel function (cf. [6, formula (1.103))

$$ K_\nu(x)K_\nu(y)= {1\over 2}\int_0^\infty e^{-{1\over 2}\left(t{x^2+y^2\over xy}+{xy\over t}\right)}K_\nu(t){dt\over t}\eqno(1.10)$$
to write the integral representation for the square of the Macdonald function $K_{i\tau}^2 (x)$ which is nonnegative for real $\tau$ and positive $x$

$$ K_{i\tau}^2 (x) = {1\over 2}\int_0^\infty e^{- t - {x^2\over 2t}}K_{i\tau}(t){dt\over t}.\eqno(1.11)$$
Hence, appealing to Entry 2.2.1.8 in [1] and (1.1) we calculate the Fourier cosine transform by $\tau$ of the left-hand side in (1.11), interchanging the order of integration on the right-hand side via the dominated convergence. Thus we obtain

$$\int_0^\infty \cos(\tau y)K^2_{i\tau}(x) d\tau = {\pi\over 2} K_0\left( 2x\cosh\left({y\over 2}\right)\right).\eqno(1.12)$$
Fulfilling again the termwise differentiation with respect to $y$, we find 

$$\int_0^\infty \cos(\tau y) \tau^{2n} K^2_{i\tau}(x) d\tau =  (-1)^n {\pi\over 2}  D_y^{2n} \bigg[ K_0\left( 2x\cosh\left({y\over 2}\right)\right)\bigg],\eqno(1.13)$$
as well as a polynomial extension

$$\int_0^\infty \cos(\tau y) q_n(\tau^{2})  K^2_{i\tau}(x) d\tau =   {\pi\over 2} \ q_n\left( - D_y^2\right) \bigg[  K_0\left( 2x\cosh\left({y\over 2}\right)\right)\bigg].\eqno(1.14)$$
In the meantime, the Mellin-Barnes representation of  $K^2_{i\tau}(x)$ yields (see [1, Entry 3.14.18.4])

$$K^2_{i\tau}(x) = {1\over 8\sqrt \pi i} \int_{\gamma-i\infty}^{\gamma+i\infty}  \Gamma\left({s\over 2} +i\tau\right) \Gamma\left({s\over 2} -i\tau\right) {\Gamma(s/2)\over \Gamma((s+1)/2) }  x^{-s} ds,\ \gamma >0.\eqno(1.15)$$
Hence owing to the absolute and uniform convergence and the addition formula for the gamma function we have

 $$ x{d\over dx} x {d\over dx} \left[ K^2_{i\tau}(x) \right]  = {1\over 8\sqrt \pi i} \int_{\gamma-i\infty}^{\gamma+i\infty}  \Gamma\left({s\over 2} +i\tau\right) \Gamma\left({s\over 2} -i\tau\right) {\Gamma(s/2)\over \Gamma((s+1)/2) }   s^2 x^{-s} ds$$

$$= {1\over 2\sqrt \pi i} \int_{\gamma-i\infty}^{\gamma+i\infty}  \Gamma\left({s\over 2} +i\tau+1\right) \Gamma\left({s\over 2} -i\tau+1\right) {\Gamma(s/2)\over \Gamma((s+1)/2) }    x^{-s} ds- \tau^2  K^2_{i\tau}(x) $$

$$= {1\over 2\sqrt \pi i} \int_{\gamma-i\infty}^{\gamma+i\infty}  \Gamma\left({s\over 2} +i\tau+1\right) \Gamma\left({s\over 2} -i\tau+1\right) {\Gamma(s/2+1)\over \Gamma((s+3)/2) }    x^{-s} ds$$

$$+  {1\over 2\sqrt \pi i} \int_{\gamma-i\infty}^{\gamma+i\infty}  \Gamma\left({s\over 2} +i\tau+1\right) \Gamma\left({s\over 2} -i\tau+1\right) {\Gamma(s/2+1)\over \Gamma((s+3)/2) }    {x^{-s}\over s} ds- \tau^2  K^2_{i\tau}(x) $$

$$= 4 x^2  K^2_{i\tau}(x) -  \tau^2  K^2_{i\tau}(x) + 4\int_x^\infty  y K^2_{i\tau}(y) dy,$$
i.e.  $K^2_{i\tau}(x)$ is the eigenfunction of the following integro-differential operator

$$\left({\cal B}  f \right) (x) \equiv  \left( 4 x^2 - x{d\over dx} x {d\over dx} \right) f(x)  + 4\int_x^\infty  y f(y) dy\eqno(1.16)$$
and we find

$$ {\cal B}  \left[K^2_{i\tau}(x)  \right] = \tau^2  K^2_{i\tau}(x).\eqno(1.17)$$
Consequently, we get from (1.13) the polynomial operator equality

$$h_n\left( {\cal B}\right) \bigg[ K_0\left( 2x\cosh\left({y\over 2}\right)\right) \bigg] =  h_n\left( - D_y^2\right) \bigg[ K_0\left( 2x\cosh\left({y\over 2}\right)\right) \bigg].\eqno(1.18)$$
Now it is worth to mention the sequence of polynomials $\{p_n(x)\}_{n\ge 0}$ related to the Kontorovich-Lebedev transform [7] which is defined in terms of the $n$th iteration of the differential operator (1.5)
$$p_n(x)= (-1)^n e^{x}{\mathcal A}^n\  e^{-x}, \ n \in \mathbb{N}_0.\eqno(1.19)$$
It has the integral representation
$$p_n(x)= {2(-1)^n\over \pi} e^{x}\int_{0}^\infty \tau^{2n} K_{i\tau}(x)\ d\tau,\ x >0,\eqno(1.20) $$
 and satisfies  the differential recurrence relation of the form
$$p_{n+1}(x)= x^2p_n^{\prime\prime}(x) + x(1-2x)p_n^\prime(x)- xp_n(x), \ n= 0,1,2,\dots \ .\eqno(1.21)$$
The generating function for this system of polynomials is given by the series
$$e^{- 2x\sinh^2( y/2)} = \sum_{n=0}^\infty {p_n(x)\over (2n)!}\ y^{2n}\eqno(1.22)$$
with a positive convergence radius.  Hence

$$e^{- x\cosh y} = e^{-x} \sum_{k=0}^\infty {p_k(x)\over (2k)!}\ y^{2k},$$
and via termwise differentiation by $y$ within the convergence radius, we derive

$$D_y^{2n} \bigg[ e^{-x\cosh y}\bigg] = e^{-x} \sum_{k=0}^\infty {p_{k+n}(x)\over (2k)!}\ y^{2k}.\eqno(1.23)$$
The latter equality proves (1.7), employing (1.19). Returning to (1.20),  we write

$$ \int_{2x}^\infty {e^{-t}\ p_n(t) \over (t^2-4x^2)^{1/2}} \ dt = {2(-1)^n\over \pi} \int_{2x}^\infty {1 \over (t^2-4x^2)^{1/2}} \int_0^\infty \tau^{2n} K_{i\tau}(t) d\tau dt.$$
The interchange of the order of integration on the right-hand side of the latter equality is allowed by the dominated convergence.  Then we appeal to Entry 2.16.3.6 in [4, Vol. II] to arrive at the formula for the moments of the square of the Macdonald function 

$$\mu_n(x)= \int_0^\infty \tau^{2n} K_{i\tau}^2 (x) d\tau = {(-1)^n \pi\over  2^{2n+1}}  \int_{1}^\infty {e^{-2xt}\ p_n(2xt) \over (t^2- 1)^{1/2}} \ dt$$

$$= {(-1)^n \pi\over  2^{2n+1}}  \int_{0}^\infty {e^{-2x\cosh u}}\ p_n(2x\cosh u) \ du,\quad x >0,\ n \in \mathbb{N}_0.\eqno(1.24)$$
Since $p_0(x)= 1$, we make a simple substitution the reciprocal Fourier cosine transform in (1.1) to get 

$$\mu_0(x)= {\pi\over 2} K_0(2x).$$
Hence from (1.22) we establish  the generating function for the moments $\mu_n$

$$ \sum_{n=0}^\infty { (-1)^n \mu_n(x)\over (2n)!}\ (2y)^{2n} = {\pi\over 2} \sum_{n=0}^\infty {y^{2n} \over (2n)!}\  \  \int_{1}^\infty {e^{-2xt}\ p_n(2xt) \over (t^2- 1)^{1/2}} \ dt $$

$$= {\pi\over 2} \int_{1}^\infty {e^{-2xt\cosh y} \over (t^2- 1)^{1/2}} \ dt =  {\pi\over 2}  K_0(2x\cosh y),$$
where the interchange of the summation and integration is possible within the interval  $|y| < \pi/4$ (cf. [7]) due to the absolute and uniform convergence.  Moreover, it confirms (1.13).  On the other hand, we employ (1.11) to derive

$$\mu_n(x)= {1\over 2} \int_0^\infty \tau^{2n} \int_0^\infty e^{- t - {x^2\over 2t}}K_{i\tau}(t){dt\over t} d\tau $$

$$= {1\over 2}  \int_0^\infty e^{- t - {x^2\over 2t}} \int_0^\infty \tau^{2n} K_{i\tau}(t) d\tau {dt\over t} $$ 

$$= {\pi (-1)^n \over 4}  \int_0^\infty e^{- 2t - {x^2\over 2t}} p_n(t) {dt\over t},\eqno(1.25)$$
where the interchange of the order of integration is owing to Fubini's theorem.  Hence,  involving  the  recurrence relation for polynomials $p_n$  [7]

$$p_{n+1}(x)= -x \sum_{k=0}^n \binom{2n+1}{2k} p_k(x),\quad n \in \mathbb{N}_0,\eqno(1.26)$$
we find

$$\mu_n(x)=  {\pi (-1)^{n+1} \over 4}  \sum_{k=0}^{n-1} \binom{2n-1}{2k}  \int_0^\infty e^{- 2t - {x^2\over 2t}} p_k(t) dt,\quad n \in \mathbb{N}.$$
Therefore the  derivative $\mu_n^\prime(x)$ has the form, accordingly,

$$\mu^\prime_n(x)= {\pi (-1)^{n} x \over 4}  \sum_{k=0}^{n-1} \binom{2n-1}{2k}  \int_0^\infty e^{- 2t - {x^2\over 2t}} {p_k(t)\over t} dt$$

$$= x \sum_{k=0}^{n-1} (-1)^{n+k} \binom{2n-1}{2k} \mu_k(x).\eqno(1.27)$$
Meanwhile, from (1.25) we have

$$\mu^\prime_n(x)= {\pi (-1)^{n+1} x \over 4}  \int_0^\infty e^{- 2t - {x^2\over 2t}} p_n(t) {dt\over t^2}.$$
Integrating by parts, we get

$${\pi (-1)^{n+1} x \over 4}  \int_0^\infty e^{- 2t - {x^2\over 2t}} p_n(t) {dt\over t^2} = {\pi (-1)^{n+1}  \over  x}  \int_0^\infty e^{- 2t - {x^2\over 2t}} p_n(t) dt $$

$$+ {\pi (-1)^{n}  \over 2x}  \int_0^\infty e^{- 2t - {x^2\over 2t}} p^\prime_n(t) dt.$$
Using the relation for derivatives  $p^\prime_n$  [7]

$$p^\prime_{n}(x)= - \sum_{k=0}^{n-1} \binom{2n}{2k} p_k(x),\quad n \in \mathbb{N},$$
it gives 

$$\mu^\prime_n(x)=  {\pi (-1)^{n+1}  \over  2x}  \bigg[ \int_0^\infty e^{- 2t - {x^2\over 2t}} p_n(t) dt +  \sum_{k=0}^{n} \binom{2n}{2k} \int_0^\infty e^{- 2t - {x^2\over 2t}} p_k(t) dt\bigg].$$
Comparing with (1.27), we deduce the equality

$$\sum_{k=0}^{n-1} (-1)^{k+1} \binom{2n-1}{2k} \mu_k(x) = {\pi  \over  2x^2} \bigg[  \int_0^\infty e^{- 2t - {x^2\over 2t}} p_n(t) dt\bigg.$$

$$\bigg.  +   \sum_{k=0}^{n} \binom{2n}{2k} \int_0^\infty e^{- 2t - {x^2\over 2t}} p_k(t) dt\bigg].$$
But from (1.25) we observe that

$${\pi (-1)^n \over 4}  \int_0^\infty e^{- 2t - {x^2\over 2t}} p_n(t) dt = \int_x^\infty  y \mu_n(y) dy,\quad n \in \mathbb{N}_0.$$
Consequently,

$$-x^2\mu_0(x) + x^2 \sum_{k=1}^{n-1} (-1)^{k+1} \binom{2n-1}{2k} \mu_k(x) = 2  \bigg[  (-1)^n \int_x^\infty  y \mu_n(y) dy\bigg.$$

$$\bigg.  +   \sum_{k=0}^{n} (-1)^k \binom{2n}{2k} \int_x^\infty  y \mu_k (y) dy\bigg],\ n \in \mathbb{N}.\eqno(1.28)$$
Differentiating through  (1.28) and using (1.26), it implies after simplification 

$$x\mu_0^\prime(x) +  x^2 \sum_{k=0}^{n-2} \binom{2n-1}{2(k+1)} \sum_{m=0}^{k} (-1)^{m} \binom{2k+1}{2m} \mu_m(x)$$

$$ = 4  (-1)^n  \mu_n(x) + 2\sum_{k=0}^{n-2} (-1)^{k+1} \binom{2n-1}{2k+1} \mu_{k+1} (x),\quad  n \in \mathbb{N}.$$
Changing the summation, it yields the following recurrence relation for the moments $\mu_n$ 

$$x\mu_0^\prime(x) +  x^2 (2n-1) \sum_{k=0}^{n-2} (-1)^{k} \binom{2(n-1)}{2k}   \mu_k(x) \sum_{m=0}^{n-k-2}  {1\over 2(k+m+1)} \binom{2(n-k-1)}{2m+1} $$

$$ = 4  (-1)^n  \mu_n(x) + 2\sum_{k=1}^{n-1} (-1)^{k} \binom{2n-1}{2k-1} \mu_{k} (x),\quad x > 0,\  n \in \mathbb{N}.\eqno(1.29)$$
Letting $n=1$, it is easily seen that $x\mu_0^\prime(x) = - 4\mu_1(x)$.  Furthermore, the inner sum on the left-hand side of (1.29) can be calculated in terms of the Gauss hypergeometric function ${}_2F_1$, employing relation 4.2.3.3 in [4, Vol. I]. Precisely, we derive

$$\sum_{m=0}^{n-k-2}  {1\over 2(k+m+1)} \binom{2(n-k-1)}{2m+1} = {1\over 2} \int_0^1 x^{2k} \bigg[ (1+x)^{2(n-k-1)} - (1-x)^{2(n-k-1)}\bigg] dx $$

$$=  {1\over 2} \int_0^1 x^{2k} (1+x)^{2(n-k-1)} dx  - {(2k)! (2(n-k-1))!\over 2 (2n-1)!}$$

$$=  {1\over 2 (2k+1) } \bigg[  {}_2F_1\left(- 2(n-k-1),\  2k+1; \  2(k+1); - 1\right)   - {(2k+1)! (2(n-k-1))!\over  (2n-1)!}\bigg],$$
where Entry 2.2.6.1 in [4, Vol. I] is applied.  Thus we rewrite (1.29) in the form

$$ (-1)^n  \mu_n(x) + \left( {3\over 2} -n\right) \mu_1(x) + {1\over 2} \sum_{k=2}^{n-1} (-1)^{k} \binom{2n-1}{2k-1} \mu_{k} (x) $$

$$ = {x^2\over 8} \sum_{k=0}^{n-2} (-1)^{k}  \bigg[ \binom{2n-1}{2k+1} {}_2F_1\left( 2(k+1-n),\  2k+1; \  2(k+1); - 1\right)   - 1\bigg] \mu_k(x),\eqno(1.30)$$
where  $x > 0,\  n \in \mathbb{N}.$

It is worth to mention the Lebedev-Skalskaya type kernels (see [6, Ch. 6]) and their Fourier cosine and sine transforms [6, formulae (6.26), (6.27)], namely,

$$\int_0^\infty \cos(\tau y) {\rm Re} \left[ K_{\alpha+i\tau}(x) \right] d\tau = {\pi\over 2}\ e^{-x\cosh y} \cosh\left(\alpha y\right),\quad x >0,\ y, \alpha \in \mathbb{R},\eqno(1.31)$$

$$\int_0^\infty \sin(\tau y) {\rm Im} \left[ K_{\alpha+i\tau}(x) \right] d\tau = {\pi\over 2}\ e^{-x\cosh y} \sinh (\alpha y),\quad x >0,\ y,\alpha \in \mathbb{R}.\eqno(1.32)$$
Hence as above

$$\int_0^\infty \cos(\tau y) \tau^{2n}\  {\rm Re} \left[ K_{\alpha+i\tau}(x) \right] d\tau = {\pi\over 2}\ (-1)^n  D^{2n}_y \bigg[e^{-x\cosh y} \cosh\left(\alpha y\right)\bigg],\eqno(1.33)$$

$$\int_0^\infty \sin(\tau y) \tau^{2n}\  {\rm Im} \left[ K_{\alpha+i\tau}(x) \right] d\tau = {\pi\over 2}\ (-1)^{n}  D^{2n}_y \bigg[ e^{-x\cosh y} \sinh (\alpha y)\bigg].\eqno(1.34)$$
Meanwhile, (1.2), (1.34) and Entries 4.2.3.2, 4.2.3.3 in [4, Vol. I]  imply the representations 

$$D^{2n}_y \bigg[e^{-x\cosh y} \cosh\left(\alpha y\right)\bigg] = \cosh\left(\alpha y\right) \sum_{k=0}^{n}\binom{2n}{2k}  D^{2k}_y \bigg[e^{-x\cosh y}\bigg]   \alpha^{2(n-k)}$$

$$+ \sinh\left(\alpha y\right) \sum_{k=0}^{n-1}\binom{2n}{2k+1}  D^{2k+1}_y \bigg[e^{-x\cosh y}\bigg]   \alpha^{2(n-k)-1}$$

$$= {2\over \pi}  \bigg[ \cosh\left(\alpha y\right) \int_0^\infty \cos(\tau y)  {\rm Re} \left[ \left(\alpha+ i\tau\right)^{2n} \right] K_{i\tau} (x) d\tau \bigg.$$

$$\bigg. -x   \sinh\left(\alpha y\right) \int_0^\infty {\sin (\tau y)\over \tau}   {\rm Im} \left[ \left(\alpha+ i\tau\right)^{2n}\right] {\rm Im} \left[ K_{1+i\tau}(x) \right] d\tau \bigg],$$

 $$D^{2n}_y \bigg[ e^{-x\cosh y} \sinh (\alpha y)\bigg] = {2\over \pi}  \bigg[ \sinh\left(\alpha y\right) \int_0^\infty \cos(\tau y)  {\rm Re} \left[ \left(\alpha+ i\tau\right)^{2n} \right] K_{i\tau} (x) d\tau \bigg.$$

$$\bigg. -x  \cosh\left(\alpha y\right) \int_0^\infty {\sin (\tau y)\over \tau}   {\rm Im} \left[ \left(\alpha+ i\tau\right)^{2n}\right] {\rm Im} \left[ K_{1+i\tau}(x) \right] d\tau \bigg].$$
But the properties of the Macdonald function [4, Vol. II] suggest the equality

$$x  {\rm Im} \left[ K_{1+i\tau}(x) \right]  = \tau K_{i\tau} (x).\eqno(1.35)$$
Consequently, we have final representations for the consecutive derivatives

 $$D^{2n}_y \bigg[e^{-x\cosh y} \cosh\left(\alpha y\right)\bigg] =   {2\over \pi}  \int_0^\infty   {\rm Re} \left[ \cosh\left( (\alpha + i\tau )y\right) \left(\alpha+ i\tau\right)^{2n} \right] K_{i\tau} (x) d\tau,\eqno(1.36)$$

$$D^{2n}_y \bigg[ e^{-x\cosh y} \sinh (\alpha y)\bigg] =  {2\over \pi} \int_0^\infty   {\rm Re} \left[ \sinh\left( (\alpha + i\tau )y\right) \left(\alpha+ i\tau\right)^{2n} \right] K_{i\tau} (x) d\tau.\eqno(1.37)$$

Further, the convolution operator for the Kontorovich-Lebedev transform is defined as follows [5]
$$(f*g)(x) \equiv f(x)* g(x) ={1\over 2x}\int_0^\infty\int_0^\infty e^{-{1\over 2}\left(x{u^2+y^2\over uy}+{yu\over x}\right)}f(u)g(y)dudy,\quad x>0.\eqno(1.38)$$
It is well defined in the Banach ring $L^\nu ({\mathbb{R}}_+) \equiv L_{1}({\mathbb{R}}_+;K_{\nu}(x)dx), \nu \in {\mathbb{R}}$, i.e. the space of all summable functions $f:
{\mathbb{R}}_+ \to {\mathbb{C}}$ with respect to the measure $K_{\nu}(x)dx$ for which
 $$||f||_{L^\nu ({\mathbb{R}}_+)}=\int_0^\infty |f(x)|K_{\nu} (x)dx$$
is finite.  Moreover,  the factorization  property is true via the product formula (1.10) in terms of the Kontorovich-Lebedev transform 

$$(Ff) (\tau)= \int_0^\infty K_{i\tau}(x) f(x) dx\eqno(1.39)$$
in the space $L^\nu({\mathbb{R}}_+)$. Namely, we have
 $$(F (f*g))(\tau)=(Ff)(\tau)(Fg)(\tau),\  \tau \in {\mathbb R}$$
when $f, g \in L^\nu({\mathbb{R}}_+)$.  In particular, for powers  $f(x)= x^{2a-1},\ g(x)= x^{2b-1},\ a,b >  |\nu|/2$  we have, employing Entries 2.3.16.1 in [4, Vol. I] and 3.14.2.2 in [1]

$$x^{2a-1}* x^{2b-1} =  {1\over 2 x} \int_0^\infty\int_0^\infty e^{-{u (x^2+y^2)\over 2 xy} - {xy\over 2u}}\  u^{2a-1} y^{2b-1} dudy$$

$$=  x^{2a-1}  \int_0^\infty   {K_{2a} \left(\sqrt{x^2+y^2}\right)\over  (x^2+y^2)^{a}} \  y^{2(a+b)-1} dy$$

$$=  x^{2(a+b)-1}\ {1\over 2}  \int_0^\infty   {K_{2a} \left(x\sqrt{1+y}\right)\over  (1+y)^{a}} \  y^{a+b-1} dy$$

$$=  2^{a+b-1}\ x^{a+b-1}\  \Gamma(a+b)  K_{b-a} (x).\eqno(1.40)$$
On the other hand, due to the Parseval equality for the Kontorovich-Lebedev transform [6, formula (4.34)] and Entry 3.14.1.3 in [1]  it yields from the previous equalities

$$ x^{a+b-1}\  \Gamma(a+b)  K_{b-a} (x) = {2^{a+b-2} \over \pi^2} \int_0^\infty \tau\sinh(\pi\tau) {K_{i\tau}(x)\over  x} \bigg| \Gamma\left( a +{i\tau\over 2}\right) \Gamma\left( b +{i\tau\over 2}\right)\bigg|^2 d\tau.$$
Hence, multiplying both sides by $x^{2c},\ c >0$ and integrating through over $\mathbb{R}_+$, we find

$$ {1 \over 4\pi^2}  \int_0^\infty \tau\sinh(\pi\tau)  \bigg| \Gamma\left( a +{i\tau\over 2}\right) \Gamma\left( b +{i\tau\over 2}\right) \Gamma\left(c +{i\tau\over 2}\right)\bigg|^2 d\tau$$

$$=  \Gamma(a+b) \Gamma(a+c) \Gamma(b+c), $$
which is  the limit case of the Wilson integral [2].  Now, considering the convolution $x^{2a-1}* x^{2b-1} * x^{2c-1},\ a,b,c > |\nu|/2$ and employing its  associativity property, we obtain

$$x^{2a-1}* x^{2b-1} * x^{2c-1} =   {2^{a+b-2}  \Gamma(a+b) \over  x} $$

$$\times \int_0^\infty\int_0^\infty e^{-{u (x^2+y^2)\over 2 xy} - {xy\over 2u}}\  u^{2c-1}   y^{a+b-1} K_{b-a} (y) dudy$$

$$=   2^{a+b-1}  \Gamma(a+b) x^{2c-1} \  \int_0^\infty\  {K_{2c} \left(\sqrt{x^2+y^2}\right)\over  (x^2+y^2)^{c}}   y^{a+b+2c-1} K_{b-a} (y) dy.$$
Meanwhile, the Parseval equality implies

$$x^{2a-1}* x^{2b-1} * x^{2c-1} $$

$$=  {2^{2(a+b+c)-5} \over \pi^2} \int_0^\infty \tau\sinh(\pi\tau) {K_{i\tau}(x)\over  x} \bigg| \Gamma\left( a +{i\tau\over 2}\right) \Gamma\left( b +{i\tau\over 2}\right)  \Gamma\left( c +{i\tau\over 2}\right)\bigg|^2 d\tau$$
and,  taking some $d > 0$, we deduce, correspondingly, 

$$\int_0^\infty \left( x^{2a-1}* x^{2b-1} * x^{2c-1}\right) x^{2d} dx $$

$$= {2^{2(a+b+c+d)-7} \over \pi^2} \int_0^\infty \tau\sinh(\pi\tau) \bigg| \Gamma\left( a +{i\tau\over 2}\right) \Gamma\left( b +{i\tau\over 2}\right)  \Gamma\left( c +{i\tau\over 2}\right) \Gamma\left( d +{i\tau\over 2}\right)\bigg|^2 d\tau.$$
But from the above and  the relation 3.14.18.4 in [1] it yields 

$$\int_0^\infty \left( x^{2a-1}* x^{2b-1} * x^{2c-1}\right) x^{2d} dx =   2^{a+b-1}  \Gamma(a+b) \int_0^\infty x^{2(c+d)-1} $$

$$\times\  \int_0^\infty\  {K_{2c} \left(\sqrt{x^2+y^2}\right)\over  (x^2+y^2)^{c}}   y^{a+b+2c-1} K_{b-a} (y) dy dx$$

$$=   2^{a+b+c+d-2}  \Gamma(a+b) \Gamma(c+d)  \int_0^\infty\   y^{a+b+c+d-1} K_{d-c} (y) K_{b-a} (y) dy $$

$$=   {2^{2(a+b+c+d)-5}  \Gamma(a+b) \Gamma(b+c) \Gamma(c+d)\Gamma(d+b) \Gamma(c+a)\Gamma(a+d) \over \Gamma (a+b+c+d)}.$$
Therefore we end up with the identity

$${1 \over 4\pi^2} \int_0^\infty \tau\sinh(\pi\tau) \bigg| \Gamma\left( a +{i\tau\over 2}\right) \Gamma\left( b +{i\tau\over 2}\right)  \Gamma\left( c +{i\tau\over 2}\right) \Gamma\left( d +{i\tau\over 2}\right)\bigg|^2 d\tau $$

$$=   { \Gamma(a+b) \Gamma(a+c) \Gamma(a+d)\Gamma(b+c) \Gamma(b+d)\Gamma(c+d) \over \Gamma (a+b+c+d)}$$
which represents the full Wilson integral [2].  Consequently, one can extend the Wilson integral for the general convolution product of powers.  Namely,  we have the formula

$$\int_0^\infty \left({}^* \prod_{k=1}^n x^{2a_k-1} \right) x^{2a_{n+1} } dx$$

$$ = {4^{\sum_{k=1}^{n+1} a_k - n} \over 2\pi^2} \int_0^\infty \tau\sinh(\pi\tau) \bigg| \prod_{k=1}^{n+1}  \Gamma\left( a_k +{i\tau\over 2}\right) \bigg|^2 d\tau,\eqno(1.41)$$
where $n= 2,3,\dots,\ a_k >  |\nu|/2,\ k =1,2,\dots, n+1.$  The right-hand side of (1.41) can be written in terms of the Meijer $G$-function [4, Vol. III] at the unity.  Indeed, we have

$$ {4^{\sum_{k=1}^{n+1} a_k - n} \over 2\pi^2} \int_0^\infty \tau\sinh(\pi\tau) \bigg| \prod_{k=1}^{n+1}  \Gamma\left( a_k +{i\tau\over 2}\right) \bigg|^2 d\tau$$

$$=  {4^{\sum_{k=1}^{n+1} a_k +1 - n} \over 2i} \int_{- i\infty}^{i\infty} { \prod_{k=1}^{n+1}  \Gamma\left( a_k + s\right)   \prod_{k=1}^{n+1}  \Gamma\left( a_k - s\right) \over \Gamma(s) \Gamma(-s)  \Gamma(1/2+s) \Gamma(1/2-s)}  ds$$

$$=  4^{\sum_{k=1}^{n+1} a_k +1 - n} \pi \mathop{G_{n+3, n+3}^{n+1, n+1}}\left( {1\bigg\vert  {1-a_1, 1-a_2,\dots, 1-a_{n+1}, 0,\  {1\over 2} 
		\atop a_1,\  a_2,\dots, \ a_{n+1}, \ 1, \ {1\over 2} }}\right).$$
Consequently,

$$\int_0^\infty \left({}^* \prod_{k=1}^n x^{2a_k-1} \right) x^{2a_{n+1} } dx=  4^{\sum_{k=1}^{n+1} a_k +1 - n} \pi $$

$$\times \mathop{G_{n+3, n+3}^{n+1, n+1}}\left( {1\bigg\vert  {1-a_1, 1-a_2,\dots, 1-a_{n+1}, 0,\  {1\over 2} 
		\atop a_1,\  a_2,\dots, \ a_{n+1}, \ 1, \ {1\over 2} }}\right).\eqno(1.42)$$
Hence, accordingly, for $y \in \mathbb{R}$

$$ {1\over 2\pi^2} \int_0^\infty \tau\sinh(\pi\tau) \cos(y\tau) \bigg| \prod_{k=1}^{n+1}  \Gamma\left( a_k +{i\tau\over 2}\right) \bigg|^2 d\tau$$

$$=   4\pi  \mathop{G_{n+3, n+3}^{n+1, n+1}}\left( {e^{-y} \bigg\vert  {1-a_1, 1-a_2,\dots, 1-a_{n+1}, 0,\  {1\over 2} 
		\atop a_1,\  a_2,\dots, \ a_{n+1}, \ 1, \ {1\over 2} }}\right).\eqno(1.43)$$
Therefore, it immediately implies the equality

$$ {1\over 2\pi^2} \int_0^\infty \tau^{2m+1} \sinh(\pi\tau)  \bigg| \prod_{k=1}^{n+1}  \Gamma\left( a_k +{i\tau\over 2}\right) \bigg|^2 d\tau$$

$$=   4\pi  (-1)^m {d^{2m} \over dy^{2m}} \bigg[ \mathop{G_{n+3, n+3}^{n+1, n+1}}\left( {e^{-y} \bigg\vert  {1-a_1, 1-a_2,\dots, 1-a_{n+1}, 0,\  {1\over 2} 
		\atop a_1,\  a_2,\dots, \ a_{n+1}, \ 1, \ {1\over 2} }}\right)\bigg] \bigg\vert_{y=0}.\eqno(1.44)$$

The latter derivative of the order $2m$ we calculate as in [8[ via Hoppe's formula.  In fact, we find, employing relation 8.2.2.32 in [4, Vol. III]

$${d^{2m} \over dy^{2m}} \bigg[ \mathop{G_{n+3, n+3}^{n+1, n+1}}\left( {e^{-y} \bigg\vert  {1-a_1, 1-a_2,\dots, 1-a_{n+1}, 0,\  {1\over 2} 
		\atop a_1,\  a_2,\dots, \ a_{n+1}, \ 1, \ {1\over 2} }}\right)\bigg] \bigg\vert_{y=0}$$

$$= \sum_{k=0}^{2m} { (-1)^k\over k!}   \mathop{G_{n+4, n+4}^{n+1, n+2}}\left( {1 \bigg\vert  {0, 1-a_1, 1-a_2,\dots, 1-a_{n+1}, 0,\  {1\over 2} 
		\atop a_1,\  a_2,\dots, \ a_{n+1}, \ 1, \ {1\over 2}, k }}\right)$$
		
$$\times  \sum_{j=0}^k (-1)^{j} \binom{k}{j}  j^{2m}.\eqno(1.45) $$		
On the other hand, recalling (1.5), (1.6) and taking  for any $j \in \{1,2,\dots, n+1\}$ the left-hand  side of (1.44) can be written in the form

$$ {1\over 2\pi^2} \int_0^\infty \tau^{2m+1} \sinh(\pi\tau)  \bigg| \prod_{k=1}^{n+1}  \Gamma\left( a_k +{i\tau\over 2}\right) \bigg|^2 d\tau$$		

$$=  {2^{1-2a_j} \over \pi^2} \int_0^\infty x^{2a_j-1}  \int_0^\infty \tau^{2m+1} \sinh(\pi\tau)  K_{i\tau}(x) \bigg| \prod_{k=1, k\neq j}^{n+1}  \Gamma\left( a_k +{i\tau\over 2}\right) \bigg|^2 d\tau dx$$		

$$=  {2^{1-2a_j} \over \pi^2} \int_0^\infty x^{2a_j-1}  \int_0^\infty \tau \sinh(\pi\tau)  {\cal A}^m K_{i\tau}(x) \bigg| \prod_{k=1, k\neq j}^{n+1}  \Gamma\left( a_k +{i\tau\over 2}\right) \bigg|^2 d\tau dx.$$		
Integration by parts in the integral by $x$ and elimination of the integrated terms, it implies

$$  {2^{1-2a_j} \over \pi^2} \int_0^\infty x^{2a_j-1}  \int_0^\infty \tau \sinh(\pi\tau)  {\cal A}^m K_{i\tau}(x) \bigg| \prod_{k=1, k\neq j}^{n+1}  \Gamma\left( a_k +{i\tau\over 2}\right) \bigg|^2 d\tau dx$$

$$=  {2^{1-2a_j} \over \pi^2} \int_0^\infty   {\cal A}^m x^{2a_j}  \int_0^\infty \tau \sinh(\pi\tau)  {K_{i\tau}(x)\over x}  \bigg| \prod_{k=1, k\neq j}^{n+1}  \Gamma\left( a_k +{i\tau\over 2}\right) \bigg|^2 d\tau dx.$$	 
Finally, the Parseval identity  for the Kontorovich-Lebedev convolution [6, formula (4.34)] leads us  to the equality

$$ {1\over 2\pi^2} \int_0^\infty \tau^{2m+1} \sinh(\pi\tau)  \bigg| \prod_{k=1}^{n+1}  \Gamma\left( a_k +{i\tau\over 2}\right) \bigg|^2 d\tau$$	

$$=  4^{n- \sum_{k=0}^{n+1} a_k}  \int_0^\infty \left( x^{2a_1-1} * x^{2a_2-1} \dots x^{2a_{j-1}-1} * x^{2a_{j+1} -1} \dots * x^{2a_{n+1}-1}\right)   {\cal A}^m x^{2a_j} dx.\eqno(1.46)$$

\section{Orthogonal polynomials for the weights $K^2_{i\tau}(x), \newline  \left[ {\rm Re} \left[ K_{1+i\tau}(x) \right]\right]^2,  \left[ {\rm Im} \left[ K_{1+i\tau}(x) \right]\right]^2$}

The aim of this section is to investigate  properties and relationships of orthonormal sequences of polynomials  $ \{P_n(\tau^2, x)\}_{n\ge 0}, \  \{Q_n(\tau^2, x)\}_{n\ge 0},  \{R_n(\tau^2, x)\}_{n\ge 0}$ with a positive parameter $x >0$, satisfying the following orthogonality relations

 $$\int_{0}^\infty  P_n(\tau^2, x) P_m(\tau^2, x)  K^2_{i\tau}(x) d\tau = \delta_{n,m},\eqno(2.1)$$
 
 $$\int_{0}^\infty  Q_n(\tau^2, x) Q_m(\tau^2, x) \left[ {\rm Im} \left[ K_{1+i\tau}(x) \right]\right]^2  d\tau = \delta_{n,m},\eqno(2.2)$$

 $$\int_{0}^\infty  R_n(\tau^2, x) R_m(\tau^2, x)   \left[ {\rm Re} \left[ K_{1+i\tau}(x) \right]\right]^2 d\tau = \delta_{n,m},\eqno(2.3)$$
where $\delta_{n,m}$ the Kronecker symbol.  Starting from (2.1), we see that  up to a normalization factor, it is equivalent to the following $n$ conditions 

$$\int_{0}^\infty  P_n(\tau^2, x)  \tau^{2m}  K^2_{i\tau}(x)   d\tau = 0,\quad m= 0,1,\dots,n-1.\eqno(2.4)$$
Writing the sequence $ \{P_n\}_{n\ge 0} $ in terms of the coefficients

$$P_n(\tau, x) = \sum_{k=0}^n a_{n,k} (x) \tau^k,\eqno(2.5)$$
we denote its leading coefficient $a_{n,n}(x) \equiv a_n(x) \neq 0$, which has the value (see (2.4)) 

$$\int_{0}^\infty  P_n(\tau^2, x)  \tau^{2n}  K^2_{i\tau}(x)   d\tau = {1\over a_n(x)}\eqno(2.6)$$
and guarantees the regularity of the polynomial sequence.   Moreover, it satisfies the 3-term recurrence relation 

$$\tau  P_n(\tau, x) = A_{n+1}(x) P_{n+1} (\tau, x)+ B_n(x) P_n(\tau, x) + A_n P_{n-1} (\tau, x),\ n \in \mathbb{N}_0\eqno(2.7)$$
where $P_{-1}(\tau, x) \equiv 0$ and 

$$A_{n+1}(x)  = {a_n(x)\over a_{n+1}(x)},\quad  B_{n}(x)  = {b_n(x)\over a_{n}(x)}- {b_{n+1}(x)\over a_{n+1}(x)},\quad b_n(x) \equiv a_{n,n-1}(x).\eqno(2.8)$$
Recalling the equality (cf. (1.13), (1.14)) for the moments (1.24)

$$\mu_k(x)= \int_0^\infty \tau^{2k}  K^2_{i\tau} (x) d\tau = {\pi\over 2} {\partial^{2k}\over \partial y^{2k}} \bigg[ K_0\left(2x\cos\left({y\over 2}\right)\right)\bigg]\bigg\vert_{y=0},$$
 it immediately implies the composition operator equality for the orthogonality (2.1)

$$P_n\left( D_y^2,x\right) P_m\left(D_y^2,x\right) \bigg[ K_0\left(2x\cos\left({y\over 2}\right)\right)\bigg]\bigg\vert_{y=0}  = {2\over \pi}  \delta_{n,m}.$$
Moreover,  orthogonality conditions (2.2) become

$$P_n\left( D_y^2,x\right) D_y^{2m}  \bigg[ K_0\left(2x\cos\left({y\over 2}\right)\right)\bigg]\bigg\vert_{y=0} = 0,\quad m= 0,1,\dots,n-1.\eqno(2.9)$$
Consequently, the three-term recurrence relation (2.7) can be written in the composition form as well. Namely, we have

$$D_y^2 P_n\left( D_y^2, x\right) \bigg[ K_0\left(2x\cos\left({y\over 2}\right)\right)\bigg]\bigg\vert_{y=0} = A_{n+1} (x) P_{n+1} \left( D_y^2,x\right) \bigg[ K_0\left(2x\cos\left({y\over 2}\right)\right)\bigg]\bigg\vert_{y=0}$$

$$+ B_n(x)  P_n\left( D_y^2,x\right) \bigg[ K_0\left(2x\cos\left({y\over 2}\right)\right)\bigg]\bigg\vert_{y=0}+ A_n (x) P_{n-1}\left( D_y^2,x\right) \bigg[ K_0\left(2x\cos\left({y\over 2}\right)\right)\bigg]\bigg\vert_{y=0}.\eqno(2.10)$$
Introducing the Hankel determinants in terms of the moments $\mu_n(x)$

$$ D_{n} (x) = \begin{vmatrix}  

\mu_0(x)  & \mu_1(x) &  \dots&  \dots&   \mu_{n}(x) \\

\mu_1(x)  & \mu_2(x) &  \dots&  \dots&   \mu_{n+1}(x) \\
 
\dots&  \dots &   \dots&   \dots&   \dots \\
 
 \vdots  &   \ddots  &  \ddots & \ddots&   \vdots\\
  
\mu_n(x)  & \mu_{n+1}(x) &  \dots&  \dots&   \mu_{2n}(x)\\

 \end{vmatrix},\ n \in \mathbb{N}_0,$$
the orthonormal sequence of polynomials $\{P_n(\tau, x)\}_{n\ge 0}$ can be written, accordingly, 

$$ P_n(\tau, x)  =  {a_n(x ) \over  D_{n-1} (x)} \begin{vmatrix}  

\mu_0(x)  & \mu_1(x) &  \dots&  \dots&   \mu_{n} (x)\\

\mu_1(x)  & \mu_2(x) &  \dots&  \dots&   \mu_{n+1} (x)\\
 
\dots&  \dots &   \dots&   \dots&   \dots \\
 
 \vdots  &   \ddots  &  \ddots & \ddots&   \vdots\\
  
\mu_{n-1} (x)  & \mu_{n}(x) &  \dots&  \dots&   \mu_{2n-1}(x)\\

1  &  \tau &  \dots&  \dots&   \tau^n\\

 \end{vmatrix},\eqno(2.11)$$
 where we set $D_{-1}(x)=1$.  Recalling (1.19) and (1.25), we write (2.11) as the following $n$-fold integral, involving polynomials $p_n$ 
 
$$ P_n(\tau, x)  =  {a_n(x ) \over  D_{n-1} (x)} \left({\pi\over 4}\right)^n \int_{\mathbb{R}^n_+}  \exp\left( - 2\sum_{j=1}^n t_j - {x^2\over 2} \sum_{j=1}^n {1\over t_j} \right) $$

$$\times  \begin{vmatrix}  

1  & - p_1(t_1) &  \dots&  \dots&   (-1)^n p_{n} (t_1)\\

- p_1(t_2)  & p_2(t_2) &  \dots&  \dots&   (-1)^{n+1} p_{n+1} (t_2)\\
 
\dots&  \dots &   \dots&   \dots&   \dots \\
 
 \vdots  &   \ddots  &  \ddots & \ddots&   \vdots\\
  
(-1)^{n-1} p_{n-1} (t_n)  & (-1)^n p_{n} (t_n) &  \dots&  \dots&   - p_{2n-1}(t_n)\\

1  &  \tau &  \dots&  \dots&   \tau^n\\

 \end{vmatrix} {dt_1dt_2\dots d t_n\over t_1t_2\dots t_n}.\eqno(2.12)$$
 In particular (see (1.24)),  we have
 
 $$P_0(\tau, x) = \left[ \sqrt{{\pi\over 2} K_0\left(2x\right)}\right]^{-1/2} $$
 and since $p_1(x)= -x,\ p_2(x)= 3x^2-x$, we derive
 
 $$\mu_1(x)= {\pi x\over 4} K_1(2x),\ \mu_2(x) = {\pi x\over 16} \left[ 3x \left(K_2(2x) - K_0(2x) \right) -  K_1(2x) \right].$$
  Therefore from (2.11)
 $$ P_{1} (\tau, x) = a_1(x) \left( \tau  - {\mu_1(x)\over \mu_0(x)}\right), $$
where, employing (2.4),  the value of $a_1(x)$ takes the form 

$$a_1(x) = \sqrt{{D_0(x)\over D_1(x)}} = \left[ \mu_2(x) - {\mu^2_1(x)\over \mu_0(x)}\right]^{-1/2}$$

$$=   {4\over \sqrt{\pi x}} \left[  3x \left(K_2(2x) - K_0(2x) \right) -  K_1(2x)   -  {2  x  K_1^2(2x) \over  K_0(2x)}\right]^{-1/2}.$$
But (see [4, Vol. II])

$$ K_2(2x) - K_0(2x) = {1\over x} \ K_1(2x).$$
Thus, finally, we obtain

$$a_1(x) =   2\ \sqrt{ { 2 K_0(2x)\over \pi x K_1(2x)}}\  \bigg[ K_0(2x)   -   x  K_1(2x) \bigg]^{-1/2}.$$
Observing, that the index integral (1.24) for the moments of the square of the Macdonald function can be treated as the Lebedev type transformation [5], we appeal to its inversion to establish the integral representation for the monomials $\tau^{2n-1},\ n \in \mathbb{N}.$. 

{\bf Lemma 1}. {\it Let $n \in \mathbb{N},\ \tau \in \mathbb{R}.$ Then the following integral representation holds valid

$$\tau^{2n-1} = - {2^{2n+1} \over \pi^2}  \sinh\left({\pi\tau\over 2}\right) \int_0^{\infty}   K_{i\tau/2}\left(y\right) \left[ I_{i\tau/2} \left(y\right) + I_{-i\tau/2} \left(y\right)\right]  \mu^{\prime}_n\left( y\right) dy,\eqno(2.13)$$
where $I_\nu(z)$ is the modified Bessel function of the first kind [$4$, Vol. II].}

\begin{proof} In fact,  taking  the Kontorovich-Lebedev integral [7]

$$ {\tau^{2n-1}\over \sinh(\pi\tau)} = {(-1)^n \over \pi} \int_0^\infty e^{-t} K_{i\tau} (t) p_n(t) {dt\over t},\quad n\in \mathbb{N},\eqno(2.14)$$ 
we employ  the  Erd{\' e}lyi-Kober integral  [5] (cf. Entry 2.16.3.1 in [4, Vol.II])

$$ \int_0^{x} { K_{i\tau} (t) \over ( x^2- t^2)^{1/2}} dt= {\pi K_{i\tau/2}\left(x/2\right) \over 4\cosh(\pi\tau/2)}  \left[ I_{i\tau/2} \left({x\over 2}\right) + I_{-i\tau/2} \left({x\over 2}\right)\right],\quad x > 0\eqno(2.15)$$
 to get  via its reciprocal inversion the following  representation for the Macdonald function

$$ K_{i\tau} (x)  = {1\over 2 \cosh(\pi\tau/2)} \  {d\over dx}   \int_0^{x} { t\  K_{i\tau/2}\left(t/ 2\right)  \over ( x^2- t^2)^{1/2}}   \left[ I_{i\tau/2} \left({t\over 2}\right) + I_{-i\tau/2} \left({t\over 2}\right)\right] dt.\eqno(2.16) $$
Hence, substituting the latter expression in (2.14) and taking into account (1.24), we integrate by parts and interchange the order of integration by the dominated convergence, which can be justified with the use of the Mellin-Barnes representation for the integral on right-hand side of (2.16) (see Entry 3.14.15.4 in [1]). Consequently, it yields  the equalities 

$$ \tau^{2n-1} = { (-1)^n  \over \pi} \sinh\left({\pi\tau\over 2}\right) \int_0^\infty e^{-t} p_n(t) $$

$$\times {d\over dt}  \int_0^{t} { y\  K_{i\tau/2}\left(y/ 2\right)  \over ( t^2- y^2)^{1/2}}   \left[ I_{i\tau/2} \left({y\over 2}\right) + I_{-i\tau/2} \left({y\over 2}\right)\right]  { dy dt\over t}$$

$$= { (-1)^{n+1}  \over \pi} \sinh\left({\pi\tau\over 2}\right) \lim_{t\to 0}\   {p_n(t) \over t e^t }  \int_0^{t} { y\  K_{i\tau/2}\left(y/ 2\right)  \over ( t^2- y^2)^{1/2}}   \left[ I_{i\tau/2} \left({y\over 2}\right) + I_{-i\tau/2} \left({y\over 2}\right)\right]  dy $$

$$+ { (-1)^{n}  \over \pi} \sinh\left({\pi\tau\over 2}\right)  \int_0^{\infty}   K_{i\tau/2}\left({y\over  2}\right) \left[ I_{i\tau/2} \left({y\over 2}\right) + I_{-i\tau/2} \left({y\over 2}\right)\right] $$

$$\times {d\over dy } \int_y^\infty   {d\over dt} \bigg[  {p_n(t) \over t e^t } \bigg]   ( t^2- y^2)^{1/2} dt dy$$

$$= { (-1)^{n+1}  \over \pi} \sinh\left({\pi\tau\over 2}\right)\bigg[  \lim_{t\to 0}\   {p_n(t) \over t e^t }  \int_0^{t} { y\  K_{i\tau/2}\left(y/ 2\right)  \over ( t^2- y^2)^{1/2}}   \left[ I_{i\tau/2} \left({y\over 2}\right) + I_{-i\tau/2} \left({y\over 2}\right)\right]  dy \bigg.$$

$$\bigg. + \int_0^{\infty}   K_{i\tau/2}\left({y\over  2}\right) \left[ I_{i\tau/2} \left({y\over 2}\right) + I_{-i\tau/2} \left({y\over 2}\right)\right]  {d\over dy } \int_y^\infty   { e^{-t} p_n(t) \over ( t^2- y^2)^{1/2} }dt  dy \bigg]$$

$$= { (-1)^{n+1}  \over \pi} \sinh\left({\pi\tau\over 2}\right)\bigg[  \lim_{t\to 0}\   {p_n(t) \over t e^t }  \int_0^{t} { y\  K_{i\tau/2}\left(y/ 2\right)  \over ( t^2- y^2)^{1/2}}   \left[ I_{i\tau/2} \left({y\over 2}\right) + I_{-i\tau/2} \left({y\over 2}\right)\right]  dy \bigg.$$

$$\bigg. + {2^{2n+1} (-1)^n \over \pi} \int_0^{\infty}   K_{i\tau/2}\left(y\right) \left[ I_{i\tau/2} \left(y\right) + I_{-i\tau/2} \left(y\right)\right]  \mu^{\prime}_n\left( y\right) dy \bigg].$$
In the meantime, recalling (2.16), we observe that

$$\lim_{t\to 0}\   {p_n(t) \over t e^t }  \int_0^{t} { y\  K_{i\tau/2}\left(y/ 2\right)  \over ( t^2- y^2)^{1/2}}   \left[ I_{i\tau/2} \left({y\over 2}\right) + I_{-i\tau/2} \left({y\over 2}\right)\right]  dy $$

$$= 2\cosh\left({\pi \tau \over 2}\right) \lim_{t\to 0}\   {p_n(t) \over t e^t }  \int_0^t K_{i\tau} (y) dy = 0,\ n \in \mathbb{N}.$$
Thus we prove equality (2.13) and complete the proof of Lemma 1. 

\end{proof}

This lemma applies to establish an  integral operator  whose eigenfunctions are arbitrary polynomials. Precisely, we have the following result.

{\bf Theorem 1}.  {\it  Let $n \in \mathbb{N},\ \tau \in \mathbb{R}.$  Any polynomial $f_n(x) = \sum_{k=0}^n f_{n,k}x^k$ satisfies the second kind integral equation of the form}

$$ f_n(4\tau^2) - f_{n,0} = {i\tau  \over 2} \ \sinh\left(\pi\tau\right)  \int_{-\infty}^{\infty} {f_n(4(y-i)^2) - f_{n,0} \over  (y-i) [\cosh(\pi y)+ \cosh(\pi\tau)]}  dy.\eqno(2.17)$$

\begin{proof} Indeed, taking into account (1.26), we write (2.13) in the form

$$\tau^{2n-1} = {2^{2n+1} (-1)^{n+1} \over \pi^2}  \sinh\left({\pi\tau\over 2}\right) \sum_{k=0}^{n-1} (-1)^{k} \binom{2n-1}{2k}  $$

$$\times \int_0^{\infty}  y  K_{i\tau/2}\left(y\right) \left[ I_{i\tau/2} \left(y\right) + I_{-i\tau/2} \left(y\right)\right]  \mu_k \left( y\right) dy.$$
Employing (1.24) and relation 4.2.3.2 in [4, Vol. I], it reads

$$\tau^{2n-1} =  {i  \over \pi^2}   \sinh\left({\pi\tau\over 2}\right) \int_0^{\infty}  y  K_{i\tau/2}\left(y\right) \left[ I_{i\tau/2} \left(y\right) + I_{-i\tau/2} \left(y\right)\right] $$

$$\times \int_{-\infty} ^{\infty}  (t-2i)^{2n-1}  K^2_{it/2} (y) dt dy.\eqno(2.18)$$
The interchange of the order of integration is guaranteed by the inequality (1.3) and the dominated convergence.  Hence we  obtain the integral equation for monomials

$$\tau^{2n} =     \int_{-\infty} ^{\infty} (t-2i)^{2n-1}   {\cal K} (\tau, t)  dt,\quad n \in \mathbb{N},\eqno(2.19)$$
where the kernel $ {\cal K}$ is defined by the integral

$$ {\cal K} (\tau, y)  =   { i \over 4\pi^2}   \tau \sinh\left({\pi\tau\over 2}\right) \int_0^{\infty}  u K_{i\tau/2}\left({u\over 2}\right) K^2_{iy/2} \left({u\over 2}\right) \left[ I_{i\tau/2} \left({u\over 2}\right) + I_{-i\tau/2} \left({u\over 2}\right)\right] du.\eqno(2.20) $$
Meanwhile, the kernel (2.20) can be represented,  recalling Entry 2.16.3.6 in [4, Vol. II]. Then with the use of (2.16) and Fubini's theorem, we derive

$$  {\cal K} (\tau, y)  =   { i \over 2\pi^2}   \tau \sinh\left({\pi\tau\over 2}\right) \int_0^{\infty}  u K_{i\tau/2}\left({u\over 2}\right) \int_u^\infty  {K_{iy} (t) \over ( t^2-u^2)^{1/2}} dt \left[ I_{i\tau/2} \left({u\over 2}\right) + I_{-i\tau/2} \left({u\over 2}\right)\right] du $$

$$= { i \over 2\pi^2}   \tau \sinh\left(\pi\tau\right) \int_0^{\infty} K_{iy} (t) \int_0^t  K_{i\tau} (u) du dt. \eqno(2.21) $$
The latter iterated integral we calculate, employing relations 2.4.18.4, 2.5.46.5, 2.5.46.9 in [4, Vol. I], 2.16.2.1, 2.16.6.1 in [4, Vol. II] under the same justifications.   Hence we derive

$$  {\cal K} (\tau, y)  = { i \over 2\pi^2}   \tau \sinh\left(\pi\tau\right) \int_0^{\infty} K_{iy} (t) \int_0^t  K_{i\tau} (u) du dt $$

$$= { i \over 2\pi^2}   \tau \sinh\left(\pi\tau\right) \int_0^{\infty} K_{iy} (t) \int_0^\infty \bigg[ 1- e^{-t\cosh u} \bigg] {\cos (\tau u) \over \cosh u} du dt $$

$$= { i \over 2\pi}   \tau \sinh\left(\pi\tau\right) \bigg[ {\pi\over 4 \cosh (\pi y/2) \cosh (\pi\tau/2)} - {2\over \sinh(\pi y)}  \int_0^{\infty} {\sin(uy)\cos (\tau u) \over \sinh (2u) } du \bigg]$$

$$=  {1\over 4} \ {i\tau \sinh\left(\pi\tau/2\right) \over   \cosh(\pi y/2)+ \cosh(\pi\tau/2)}.\eqno(2.22)$$
Thus  we easily end up with (2.17) after simple substitutions,  completing the proof of Theorem 1. 

\end{proof}

{\bf Corollary 1}. {\it Let $\tau, y \in \mathbb{R}.$ One has the value of the following  integral with the product of the modified Bessel functions

$$   \int_0^\infty  x K_{i\tau}\left(x\right) K^2_{iy} \left(x\right) \left[ I_{i\tau} \left(x\right) + I_{-i\tau} \left(x\right)\right] dx $$

$$= \int_0^\infty  x K_{iy}\left(x\right) K^2_{i\tau} \left(x\right) \left[ I_{iy} \left(x\right) + I_{-iy} \left(x\right)\right] dx =  {\pi^2 \over  4[ \cosh(\pi y)+ \cosh(\pi\tau)]}.\eqno(2.23) $$
Reciprocally, via the Lebedev transform (cf. [5, formula (7.53)]), the square of the Macdonald function  can be represented by the index integral}

$$ x K^2_{i\tau} \left(x\right) = - {d\over dx}  \int_0^\infty  {y\sinh(\pi y) K^2_{i y} (x) \over \cosh(\pi \tau)+ \cosh(\pi y) } d y,\ x >0, \tau \in \mathbb{R}.\eqno(2.24)$$
Writing (2.17) for polynomials (2.3), it reads

$$ P_n(4\tau^2, x) - a_{n,0} (x) = {i\tau  \over 2} \ \sinh\left(\pi\tau\right) \int_{-\infty}^{\infty}   { P_n(4(y-i)^2,x) - a_{n,0}(x)\over (y-i)[\cosh(\pi y)+ \cosh(\pi\tau)]} dy,\eqno(2.25)$$
where $\ x >0, \ \tau \in \mathbb{R}.$ Then the 3-term recurrence relation (2.5) implies

$$ P_n(4\tau^2, x)  = {i \over 2\tau} \ \sinh\left(\pi\tau\right) \int_{-\infty}^{\infty}   { (y-i) P_n(4(y-i)^2,x) \over \cosh(\pi y)+ \cosh(\pi\tau)} dy.\eqno(2.26)$$
Meanwhile, since the moments (1.24) are continuously differentiable functions for $ x >0$,  and via (2.11) we conclude the same property for polynomial coefficients $a_{n,k}(x)$ (see (2.3)),
one returns to orthogonality conditions (2.2) to differentiate this equality with respect to $x$. Hence it yields

 $$\int_{0}^\infty  {\partial P_n(\tau^2, x)\over \partial x}   \tau^{2m}   K^2_{i \tau} (x)  d\tau $$
 
 $$+  2 \int_{0}^\infty  P_n(\tau^2, x)  \tau^{2m}   K_{i \tau} (x) {\partial  K_{i \tau} (x) \over \partial x} d\tau= 0,\ m=0,1,\dots,n-1,\eqno(2.27)$$
where the differentiation under the integral sign can be justified, employing the identity (1.35) and its companion  for the derivative of the Macdonald function  [4, Vol. II]

$$ {\partial  K_{i \tau} (x) \over \partial x} = - {\rm Re} \left[ K_{1+i\tau}(x) \right].\eqno(2.28)$$
Moreover, identities (1.35), (2.28) and simple substitutions imply the orthogonality conditions in the form

$$\int_{-\infty}^\infty  {\partial P_n(\tau^2, x)\over \partial x}   \tau^{2m}   K^2_{i \tau} (x)  d\tau  +  i x a_{n,0}(x) \int_{-\infty-i}^{\infty-i}  (\tau+i)^{2m-1}   K^2_{i \tau} (x) d\tau$$

$$ +  i x \int_{-\infty-i}^{\infty-i}  \left[ P_n((\tau+i)^2, x) -a_{n,0}(x) \right] (\tau+i)^{2m-1}   K^2_{i \tau} (x) d\tau= 0,\ m=0,1,\dots,n-1.$$
We observe that  in the latter integral it is possible  to shift the contour to the real axis due to the analyticity of the integrand in the horizontal strip $-1 \le {\rm Im} \tau \le 0$ as a function of $\tau$  and integrability properties.  In the meantime (see (1.35), (2.28), (1.24))

$$ i x a_{n,0}(x) \int_{-\infty-i}^{\infty-i}  (\tau+i)^{2m-1}   K^2_{i \tau} (x) d\tau =  i x a_{n,0}(x) \int_{0}^{\infty}  \tau^{2m-1}  \left[ K^2_{1+i \tau} (x) - K^2_{1-i \tau} (x) \right] d\tau$$

$$= 2 a_{n,0}(x) \int_{0}^{\infty}  \tau^{2m}  {\partial  K^2_{i \tau} (x) \over \partial x} d\tau = 2 a_{n,0}(x) \mu_m^\prime (x).\eqno(2.29)$$ 
The result is the same if we shift the contour to the real line in (2.29), understanding its convergence in the principal value sense when $m=0$.  Thus the orthogonality conditions (2.2) take the form 

 $$\int_{-\infty}^\infty \bigg[  {\partial P_n(\tau^2, x)\over \partial x}   \tau^{2m} +  i x   P_n((\tau+i)^2, x) (\tau+i)^{2m-1} \bigg]   K^2_{i \tau} (x)  d\tau  = 0,\ m < n.\eqno(2.30)$$
On the other hand, conditions (2.2) can be written, employing identity (1.35).  We have for $m=1,2,\dots, n-1$

$$0= \int_{-\infty}^\infty  P_n(\tau^2, x)  \tau^{2(m-1)}  \left[ \tau K_{i\tau}(x)\right]^2   d\tau =  x^2 \int_{-\infty}^\infty  P_n(\tau^2, x)  \tau^{2(m-1)}  \bigg[ {\rm Im} \left[ K_{1+i\tau}(x) \right]\bigg]^2d\tau,$$
which  means the conditions
$$\int_{-\infty}^\infty  P_n(\tau^2, x)  \tau^{2m}  \bigg[ {\rm Im} \left[ K_{1+i\tau}(x) \right]\bigg]^2  d\tau = 0,\quad m= 0,\dots,n-2.\eqno(2.31)$$
The condition (2.31) means that the sequence $\{P_n\}$ is quasi-orthogonal with respect to the weight $ \left[{\rm Im} \left[ K_{1+i\tau}(x) \right]\right]^2$. Hence, representing it in terms of $\{Q_n\}$ (see (2.2)) with connection coefficients $h_{n,k} (x)$

$$P_n(\tau, x) = \sum_{k=0}^n h_{n,k} (x) Q_k(\tau, x),$$
it has the equalities from (2.31)

$$ \sum_{k=0}^m h_{n,k} (x) d_{k,m} (x)= 0, \quad m= 0,\dots,n-2,\eqno(2.32)$$
where 
$$d_{k,m} (x)= \int_{-\infty}^\infty  Q_k(\tau^2, x)  \tau^{2m}  \bigg[ {\rm Im} \left[ K_{1+i\tau}(x) \right]\bigg]^2  d\tau.\eqno(2.33)$$
Recalling (2.6), (1.35), we find the value of the integral 

$$\int_{-\infty}^\infty  P_n(\tau^2, x)  \tau^{2(n-1)}  \bigg[ {\rm Im} \left[ K_{1+i\tau}(x) \right]\bigg]^2  d\tau = {2 \over x^2 a_n(x)}.$$
Moreover, the linear system (2.32) with a lower triangular matrix and nonzero determinant $\prod_{k=0}^m d_{k,k}(x)$,  where (see (2.33))

$$d_{k,k}(x)= \int_{-\infty}^\infty  Q_k(\tau^2, x)  \tau^{2k}  \bigg[ {\rm Im} \left[ K_{1+i\tau}(x) \right]\bigg]^2  d\tau = {2\over q_k(x)}\eqno(2.34)$$
and $q_k(x)$ is a leading coefficient of the polynomial 

$$Q_k(\tau, x) = q_k(x) \tau^k+ r_k(x) \tau^{k-1} + \hbox{lower terms},$$
has zero solution. Therefore, denoting by $\alpha_n(x)\equiv h_{n,n}(x),  \beta_n(x)\equiv h_{n,n-1}(x)$, we find the relationship between polynomials $P_n, Q_n$

$$P_n(\tau, x) =\alpha_n(x)  Q_n(\tau, x)+ \beta_n(x) Q_{n-1}(\tau, x).\eqno(2.35)$$
It is easily seen, comparing the coefficients that 

$$\alpha_n(x) = {a_n(x)\over q_n(x)},\quad    \beta_n(x) = {1\over  q_{n-1} (x)}  \left[ b_n(x)  - {a_n(x) r_n(x)\over q_n(x) }\right]. $$
Meanwhile,  (2.1), (2.2), (2.7), (1.35), (2.35)  imply

$$\int_{-\infty}^\infty  P^2_n(\tau^2, x)   \bigg[ {\rm Im} \left[ K_{1+i\tau}(x) \right]\bigg]^2  d\tau = 2[\alpha^2_n(x)+ \beta^2_n(x)]$$

$$= {1\over x^2} \int_{-\infty}^\infty  P^2_n(\tau^2, x)   \tau^2  K^2_{i\tau}(x) d\tau =  {2\over x^2} \ B_n(x),$$
Moreover, equality (2.6), (2.7), (2.34), (2.35)  yield

$$  {2 x^2 \over q_n(x)} \left[  \alpha_n(x)-  { \beta_n(x) r_n(x)\over q_{n-1} (x)} \right] = \int_{-\infty}^\infty  P_n(\tau^2, x)  \tau^{2(n+1)}  K^2_{i\tau}(x) d\tau= - { 2 b_{n+1}(x)\over a_n(x) a_{n+1}(x)},$$

$$2\alpha_n(x) = \int_{-\infty}^\infty  P_n(\tau^2, x)    Q_n(\tau^2, x)  \bigg[ {\rm Im} \left[ K_{1+i\tau}(x) \right]\bigg]^2  d\tau$$

$$ = {1\over x^2} \int_{-\infty}^\infty  P_n(\tau^2, x)  \tau^2   Q_n(\tau^2, x)   K^2_{i\tau}(x) d\tau = {2\over x^2 a_n(x) } \left[ r_n(x)  - {q_n(x) b_{n+1}(x)\over a_{n+1}(x) }\right],$$

$$2\beta_n(x) = \int_{-\infty}^\infty  P_n(\tau^2, x)    Q_{n-1} (\tau^2, x)  \bigg[ {\rm Im} \left[ K_{1+i\tau}(x) \right]\bigg]^2  d\tau$$

$$ = {1\over x^2} \int_{-\infty}^\infty  P_n(\tau^2, x)  \tau^2   Q_{n-1} (\tau^2, x)   K^2_{i\tau}(x) d\tau =  2  {q_{n-1}(x)  \over  x^2 a_{n}(x)},$$

$$ \int_{-\infty}^\infty  P_n(\tau^2, x)    Q_n(\tau^2, x)   K^2_{i\tau}(x) d\tau = {2 q_n(x)\over a_n(x)} $$

$$=  x^2 \int_{-\infty}^\infty  P_n(\tau^2, x)    {Q_n(\tau^2, x)- Q_n(0, x)\over \tau^2}  \bigg[ {\rm Im} \left[ K_{1+i\tau}(x) \right]\bigg]^2  d\tau= {2 x^2 \beta_n(x) q_n(x)\over q_{n-1}(x)}.$$ 

Thus we derived the identities

$$ B_n(x) = x^2 [\alpha^2_n(x)+ \beta^2_n(x)]=  {a^2_n(x) x^2\over q^2_n(x)}+ {q^2_{n-1}(x)  \over  x^2 a^2_{n}(x)},$$

$$  {  b_{n+1}(x)\over  a_{n+1}(x)}=   {r_n(x)\over q_n(x)} - {x^2 a_n^2(x) \over q^2_n(x)},$$

$${  b_{n}(x)\over  a_{n}(x)} = {r_n(x)\over q_n(x)}  +  {q^2_{n-1}(x)  \over  x^2 a^2_{n}(x)}.$$
Then, writing the three-term recurrence relation for polynomials $Q_n$ as follows

$$\tau Q_n(\tau, x) = \hat{A}_{n+1} (x) Q_{n+1}(\tau, x) +  \hat{B}_{n} (x) Q_{n}(\tau, x) +  \hat{A}_{n} (x) Q_{n-1}(\tau, x),$$
where $Q_{-1}(\tau, x) \equiv 0$ and 
$$\hat{A}_{n+1}(x)  = {q_n(x)\over q_{n+1}(x)},\quad  \hat{B}_{n}(x)  = {r_n(x)\over q_{n}(x)}- {r_{n+1}(x)\over q_{n+1}(x)},$$
we establish the following equality

$$ B_n(x) - \hat{B}_{n} (x) = {q^2_{n-1}(x)  \over  x^2 a^2_{n}(x)}  - {q^2_{n}(x)  \over  x^2 a^2_{n+1}(x)}= {q^2_{n}(x)  \over  x^2 a^2_{n}(x)} \left[ \hat{A}^2_{n}(x) -  A^2_{n+1}(x)\right].$$
On the other hand,  conditions (2.31) yield 
$$\int_{-\infty}^\infty  P_n(\tau^2, x)  \tau^{2m}  K^2_{1+i\tau}(x) d\tau = \int_{-\infty}^\infty  P_n(\tau^2, x)  \tau^{2m}\  {\rm Re} \left[ K^2_{1+i\tau}(x) \right] d\tau $$

$$= \int_{-\infty}^\infty  P_n(\tau^2, x)  \tau^{2m} \bigg[ {\rm Re} \left[ K_{1+i\tau}(x) \right]\bigg]^2  d\tau,\quad m= 0,\dots,n-2.$$
Furthermore,

$$\int_{-\infty}^\infty  P_n(\tau^2, x)  \tau^{2m}  K^2_{1+i\tau}(x) d\tau = \int_{-\infty+i}^{\infty+i}  P_n((\tau+i)^2, x)  (\tau+i)^{2m}  K^2_{i\tau}(x) d\tau$$

$$=  \int_{-\infty}^{\infty}  P_n((\tau+i)^2, x)  (\tau+i)^{2m}  K^2_{i\tau}(x) d\tau,$$
i.e.

$$\int_{-\infty}^{\infty}  P_n((\tau+i)^2, x)  (\tau+i)^{2m}  K^2_{i\tau}(x) d\tau $$

$$= \int_{-\infty}^\infty  P_n(\tau^2, x)  \tau^{2m} \bigg[ {\rm Re} \left[ K_{1+i\tau}(x) \right]\bigg]^2  d\tau,\quad m= 0,\dots,n-2.\eqno(2.36)$$

 {\bf Remark 1}.  The author suggests the interested readers to involve orthogonal polynomials $R_n(\tau^2, x)$ (see (2.3)) in order to establish further properties and relationships for the sequence $P_n(\tau^2, x)$.
 
 \section{Generalized Wilson polynomials }
 
 In this section, involving results of Section 1, we generalize Wilson polynomials, considering the following orthogonality
 
 $$  \int_0^\infty   \bigg|\frac{ \prod_{k=1}^{p+1}  \Gamma\left( a_k + i\tau\right)}{\Gamma(2i\tau)} \bigg|^2  W_n(\tau^2, [a]_{p+1}) W_m(\tau^2, [a]_{p+1}) d\tau = \delta_{n,m},\eqno(3.1)$$	
where $[a]_p = (a_1, a_2, \dots, a_{p+1}),\ p \in \mathbb{N}_0\backslash\{0,1\}, a_j > |\nu|/2,\ j=1,2,\dots, p+1,  \nu \in \mathbb{R}$.   The moments of the weight function are  calculated by formula (1.44).  Moreover, we find the constant polynomial  $W_0(\tau, [a]_{p+1})$ as

$$ W_0(\tau, [a]_{p+1}) = {1\over 2\sqrt \pi} \bigg[  \mathop{G_{p+3, p+3}^{p+1, p+1}}\left( {1 \bigg\vert  {1-a_1, 1-a_2,\dots, 1-a_{p+1}, 0,\  {1\over 2} 
		\atop a_1,\  a_2,\dots, \ a_{p+1}, \ 1, \ {1\over 2} }}\right)\bigg]^{-1/2}.\eqno(3.2)$$
The orthogonality conditions 
 
$$  \int_0^\infty   \bigg|\frac{ \prod_{k=1}^{p+1}  \Gamma\left( a_k + i\tau\right)}{\Gamma(2i\tau)} \bigg|^2  W_n(\tau^2, [a]_{p+1}) \tau^{2m} d\tau = 0,\ \quad m=0,1,\dots, n-1\eqno(3.3)$$
can be written via (1.46) in the form of the composition orthogonality [10], employing differential operator (1.5) and convolution (1.38).  Indeed, it becomes 

$$  \int_0^\infty \left( x^{2a_1-1} * x^{2a_2-1} \dots x^{2a_{j-1}-1} * x^{2a_{j+1} -1} \dots * x^{2a_{p+1}-1}\right) $$

$$\times  W_n\left({ {\cal A}\over 4} , [a]_{p+1} \right)  {\cal A}^m x^{2a_j} dx = 0,\quad m=0,1,\dots, n-1.\eqno(3.4)$$
However,  integrating by parts and eliminating integrated terms due to the structure of the operator (1.5), we rewrite (3.4) as follows

$$  \int_0^\infty {\cal A}^m   \bigg[ x \left( x^{2a_1-1} * x^{2a_2-1} \dots x^{2a_{j-1}-1} * x^{2a_{j+1} -1} \dots * x^{2a_{p+1}-1}\right) \bigg] $$

$$\times     W_n\left({ {\cal A} \over 4}, [a]_{p+1} \right) x^{2a_j}  {dx\over x} = 0,\quad m=0,1,\dots, n-1.\eqno(3.5)$$
 But since ${\cal A}\  x^{2a_j} = x^{2a_j} \left( x^2 - 4 a_j^2\right)$, we find
 
$$W_n\left({{\cal A}\over 4} , [a]_{p+1} \right) x^{2a_j} = x^{2a_j} \hat{W}_n (x^2 , [a]_{p+1}).\eqno(3.6) $$
Hence equalities (3.5) take the form

 $$  \int_0^\infty  x^{2a_j-1} {\cal A}^m   \bigg[ x \left( x^{2a_1-1} * x^{2a_2-1} \dots x^{2a_{j-1}-1} * x^{2a_{j+1} -1} \dots * x^{2a_{p+1}-1}\right) \bigg] $$

$$\times  \hat{W}_n (x^2 , [a]_{p+1}) \   dx = 0,\quad m=0,1,\dots, n-1.\eqno(3.7)$$
The relation between polynomials $W_n$ and $\hat{W}_n$ can be established in terms of the Kontorovich-Lebedev transform (1.39). Indeed, recalling the Parseval equality for the convolution (1.38) [6, formula (4.34)] and differentiating under the integral sign due to the absolute and uniform convergence, we get

$$ {\cal A}^m   \bigg[ x \left( x^{2a_1-1} * x^{2a_2-1} \dots x^{2a_{j-1}-1} * x^{2a_{j+1} -1} \dots * x^{2a_{p+1}-1}\right) \bigg] $$

$$= 4^{ 1/2+m-p+ \sum_{k=0, k\neq j}^{p+1} a_k}  \int_0^\infty \tau^{2m} \ K_{2i\tau} (x) \bigg| {1\over \Gamma(2i\tau)}  \prod_{k=1,\ k\neq j}^{p+1}  \Gamma\left( a_k + i\tau\right) \bigg|^2 d\tau.\eqno(3.8)$$	
Multiplying both sides of (3.8) by $x^{2a_j-1} \hat{W}_n (x^2 , [a]_{p+1})$ and integrating over $\mathbb{R}_+$, we interchange the order of integration on the right-hand of the obtained equality,  owing to the dominated convergence theorem  with the use of (1.3).   Then, taking into account (3.6),  we deduce 

$$   \int_0^\infty x^{2a_j-1} \hat{W}_n (x^2 , [a]_{p+1}) \int_0^\infty \tau^{2m} \ K_{2i\tau} (x) \bigg| {1\over \Gamma(2i\tau)}  \prod_{k=1,\ k\neq j}^{p+1}  \Gamma\left( a_k + i\tau\right) \bigg|^2 d\tau dx $$

$$= \int_0^\infty \tau^{2m} \   \int_0^\infty x^{2a_j-1} \hat{W}_n (x^2 , [a]_{p+1}) K_{2i\tau} (x) dx\  \bigg| {1\over \Gamma(2i\tau)} \bigg.$$

$$\times  \prod_{k=1,\ k\neq j}^{p+1}  \Gamma\left( a_k + i\tau\right) \bigg|^2 d\tau = 0,\quad   m=0,1,\dots, n-1.$$
Meanwhile,  employing Entry 3.14.1.3 in [1], we observe that the Kontorovich-Lebedev  integral 

$$ {1\over | \Gamma ( a_j + i\tau)|^2} \int_0^\infty x^{2a_j-1} \hat{W}_n (x^2 , [a]_{p+1}) K_{2i\tau} (x) dx\eqno(3.9) $$
is a polynomial.  Moreover, it can be calculated explicitly, employing (1.6), (3.6) and the integration by parts.  In fact, writing $W_n$ in terms of coefficients

$$W_n (\tau , [a]_{p+1}) =\sum_{k=0}^n f_{n,k} \tau^{k},\eqno(3.10)$$
we have

$$ {1\over | \Gamma ( a_j + i\tau)|^2} \int_0^\infty x^{2a_j-1} \hat{W}_n (x^2 , [a]_{p+1}) K_{2i\tau} (x) dx $$

$$= {1\over | \Gamma ( a_j + i\tau)|^2} \int_0^\infty K_{2i\tau} (x) W_n\left({{\cal A}\over 4} , [a]_{p+1} \right) x^{2a_j}    {dx\over x}  $$

$$= {1\over | \Gamma ( a_j + i\tau)|^2}  \sum_{k=0}^n f_{n,k} \int_0^\infty \tau^{2k} K_{2i\tau} (x) x^{2a_j-1}  dx = 4^{a_j-1} W_n (  \tau^2, [a]_{p+1}).\eqno(3.11)$$

Another approach is to return to equalities (3.3) and to employ the Kontorovich-Lebedev integral (2.14), substituting it on the right-hand side of (3.3). Hence we have for $m=1,2,\dots, n-1$ 

$$  \int_0^\infty \tau \sinh(\pi\tau)  \bigg|\frac{ \prod_{k=1}^{r+1}  \Gamma\left( a_k + i\tau\right)}{\Gamma(2i\tau)} \bigg|^2  W_n(\tau^2, [a]_{r+1})  \int_0^\infty e^{-x} K_{i\tau} (x) p_m(x) {dx  d\tau\over x} = 0.$$
The interchange of the order of integration is allowed for $r= 2,3,\dots$ by virtue of the dominated convergence theorem, and we write

$$  \int_0^\infty e^{-x}  p_m(x)  F_n(x) {dx \over x} = 0,\  m=1,2,\dots, n-1,\eqno(3.12)$$
where we denote by 
$$ F_n(x)=  \int_0^\infty \tau \sinh(\pi\tau)  \bigg|\frac{ \prod_{k=1}^{r+1}  \Gamma\left( a_k + i\tau\right)}{\Gamma(2i\tau)} \bigg|^2    K_{i\tau} (x) W_n(\tau^2, [a]_{r+1}) d\tau.\eqno(3.13)$$
In the meantime, appealing to inequality (1.3) it is not difficult to justify the expansion of functions $F_n(x)$ for each $n$  in series  of Laguerre polynomials [3, Ch. 4, Th. 3] to obtain

$$F_n(x) = \sum_{k=0}^\infty  c_{n,k} L_k(x),\eqno(3.14)$$
where 

$$c_{n,k} = \int_0^\infty  e^{-x} L_k(x) F_n(x) dx,\eqno(3.15)$$
and the Parseval equality holds
$$ \sum_{k=0}^\infty  c^2_{n,k} =  \int_0^\infty  e^{-x}  F^2_n(x) dx.\eqno(3.16)$$
Starting from partial sums of the series (3.14), we substitute it  into (3.12) to find via orthogonality of Laguerre polynomials,  (1.26) and the passage to limit under the integral sign (similar as in the proof of the expansion theorem) the equalities

$$  \int_0^\infty e^{-x}  p_m(x)  F_n(x) {dx \over x} =  \int_0^\infty e^{-x}  p_m(x)  \sum_{k=0}^{m-1} c_{n,k} L_k(x) {dx \over x}$$

$$= \int_0^\infty e^{-x}  x^{m-1}  \sum_{k=0}^{m-1} c_{n,k} L_k(x) dx = 0,\  m=1,2,\dots, n-1.$$
This means immediately that

$$c_{n,0}= c_{n,1} = \dots = c_{n,n-2} =0,\quad n =2,3, \dots .$$
Hence (3.14) reads as 

$$F_n(x) = \sum_{k=n-1}^\infty  c_{n,k} L_k(x)$$
and in the same manner for $m \ge n$

$$ \int_0^\infty e^{-x}  p_m(x)  F_n(x) {dx \over x} =  \int_0^\infty e^{-x}  p_m(x)  \sum_{k=n-1}^{m-1}  c_{n,k} L_k(x)  {dx \over x} $$

$$=  \sum_{k=n-1}^{m-1}  c_{n,k}  \int_0^\infty e^{-x}  p_m(x)  L_k(x)  {dx \over x} = \sum_{k=n-1}^{m-1}  c_{n,k} l_{m,k},\eqno(3.17)$$
where the known values $l_{m,k}$ are defined by
$$l_{m,k} =  \int_0^\infty e^{-x}  p_m(x)  L_k(x)  {dx \over x},\ l_{m,k} = 0,\ m < k+1,\ m,k \in \mathbb{N}.\eqno(3.18)$$
Besides, recalling (1.26), we derive, in turn,

$$ \sum_{k=n-1}^{m-1}  c_{n,k}  \int_0^\infty e^{-x}  p_m(x)  L_k(x)  {dx \over x} =  - \sum_{k=n-1}^{m-1}  c_{n,k} $$

$$\times   \sum_{q= k}^{m-1} \binom{2m-1}{2q}      \int_0^\infty e^{-x}  p_q(x)  L_k(x) dx,$$
which means 
  $$l_{m,k}= - \sum_{q= k}^{m-1} \binom{2m-1}{2q}      \int_0^\infty e^{-x}  p_q(x)  L_k(x) dx,\quad m, k \in \mathbb{N}.$$
Moreover, from  (3.18) and the 3-term recurrence relation for Laguerre polynomials [3] we find  

$$(k+1) l_{m,k+1}=  (2k+1) l_{m,k} - k l_{m,k-1} - \int_0^\infty e^{-x}  p_m(x)  L_k(x)  dx.\eqno(3.19) $$ 
Then the integration by parts and equalities for derivatives of the involved polynomials [3], [7]

$$p_n^\prime(x)= - \sum_{r=0}^{n-1}  \binom{2n}{2r} p_r(x),\eqno(3.20)$$

$$xL_n^\prime (x)= n L_n (x) - n L_{n-1} (x)$$
yield $m, k \in \mathbb{N}$

$$ \int_0^\infty e^{-x}  p_m(x)  L_k(x)  dx =  \int_0^\infty e^{-x}  p^\prime_m(x)  L_k(x)  dx + \int_0^\infty e^{-x}  p_m(x)  L^\prime_k(x)  dx $$

$$= k l_{m,k} - k l_{m,k-1}  -  \sum_{r=k-1}^{m-1}  \binom{2m}{2r} \bigg[ (k+1) l_{r,k+1}+  (2k+1) l_{r,k} + k l_{r,k-1}\bigg].$$
Thus we derive, finally, the following recurrence relation

$$l_{m,k+1} -  l_{m,k}  =   \sum_{r=k-1}^{m-1}  \binom{2m}{2r} \bigg[ l_{r,k+1}+  {2k+1\over k+1} l_{r,k} + {k\over k+1} l_{r,k-1}\bigg].$$
Further, combining with (3.17), it implies the equality  

$$ \int_0^\infty e^{-x}  p_m(x)  F_n(x) {dx \over x} =   \sum_{k=n-1}^{m-1}  c_{n,k} l_{m,k}.\eqno(3.21) $$
which is valid for all $m,n \in \mathbb{N}$, taking into account empty sums.  For $m\in \mathbb{N}$ we plug in  the right-hand side of (3.13) on the left-hand side of (3.21) and  interchange the order of integration  due  to the dominated convergence.  Then the inner integral is calculated  by [7, formula (1.20)] to obtain 

$$ \pi (-1)^m  \int_0^\infty  \bigg|\frac{ \prod_{k=1}^{r+1}  \Gamma\left( a_k + i\tau\right)}{\Gamma(2i\tau)} \bigg|^2    W_n(\tau^2, [a]_{r+1})\ \tau^{2m} d\tau= \sum_{k=n-1}^{m-1}  c_{n,k} l_{m,k}.\eqno(3.22)$$
Taking into account (3.10) and denoting by $f_{n,n} \equiv f_n,\ f_{n,n-1} \equiv \hat{f_n}$, we let $m=n$ in (3.22) and recall (3.1) to deduce

$$c_{n,n-1} = {(-1)^n \pi\over f_n\  l_{n,n-1} }.$$
Moreover,  the leading term of $p_n$ is $(- 1)^n  (2n-1)!!$ (see [7]). Hence from (3.18) we find $l_{n,n-1} =  - (2n-1) !!  (n-1)! $ and therefore

$$c_{n,n-1} = {(-1)^{n-1} \pi\over f_n\  (2n-1) !!  (n-1)!  }.\eqno(3.23)$$
Letting $m=n+1$ and since from (3.1) we have 
$$ \int_0^\infty  \bigg|\frac{ \prod_{k=1}^{r+1}  \Gamma\left( a_k + i\tau\right)}{\Gamma(2i\tau)} \bigg|^2    W_n(\tau^2, [a]_{r+1})\ \tau^{2(n+1)} d\tau= - {\hat{f}_{n+1}\over f_{n+1} f_n},$$
we find from (3.22), (3.23)
$$ {(-1)^n  \over \pi}  c_{n,n}  (2n+1) !! \ n! \ f_n+  {l_{n+1,n-1} \over (2n-1) !!  (n-1)! } = - {\hat{f}_{n+1}\over f_{n+1} }.\eqno(3.24)$$
In the meantime,  recalling (3.13), we apply operator (1.5) to both sides and differentiate under the integral sign owing to the absolute and uniform convergence.  Then (1.6) and the   3-term recurrence relation for the sequence $\{ W_n(\tau, [a]_{r+1})\}_{n\ge 0}$ written in the form 
$$ \tau W_n(\tau, [a]_{p+1}) = H_{n+1} W_{n+1} (\tau, [a]_{p+1}) + E_n W_{n} (\tau, [a]_{p+1}) + H_{n} W_{n-1} (\tau, [a]_{p+1})\eqno(3.25)$$ 
yield the following differential recurrence relation for the sequence  $\{ F_n(x)\}_{n\ge 0}$
$$ {\cal A} F_n(x)=  H_{n+1} F_{n+1} (x)  + E_n F_{n} (x) + H_{n} F_{n-1} (x),\ x > 0.\eqno(3.26)$$
Then appealing to (1.19) and making integration by parts, we  write the left-hand side of (3.21) in the form

$$ \int_0^\infty e^{-x}  p_m(x)  F_n(x) {dx \over x} = - \int_0^\infty e^{-x}  p_{m-1} (x)  {\cal A} F_n(x) {dx \over x} $$

$$= -  \int_0^\infty e^{-x}  p_{m-1} (x)\bigg[ H_{n+1}   F_{n+1}(x) + E_n F_{n} (x) + H_{n} F_{n-1} (x) \bigg]  {dx \over x}.\eqno(3.27)$$
So, comparing with the right-hand side of (3.21), the case $m=n$ immediately ensures the identity
$$ H_n =-  {l_{n,n-1} c_{n,n-1} \over  l_{n-1,n-2}\ c_{n-1,n-2}},\eqno(3.28)$$
which can be written in the form  

$$ H_n =  -  (2n-1) (n-1)\  {c_{n,n-1} \over  c_{n-1,n-2}},\ n \neq 1.\eqno(3.29)$$
Hence  via  (3.10)

$$ c_{n,n-1} = (- 1)^{n-1} {  2^{n-1} c_{1,0}  \over (2n-1) !}  \prod_{k=2}^n H_k = (- 1)^{n-1} { 2^{n-1}   c_{1,0}\ f_1\over  f_n (2n-1) !},\ H_{n+1} = {f_n\over f_{n+1}}.\eqno(3.30) $$
Comparing with (3.23),  this gives the value $  c_{1,0} f_1 = \pi$.

Another approach is to recall (3.18) and do similar to (3.27). Precisely, we have

$$l_{m,k} = -  \int_0^\infty e^{-x}  p_{m-1} (x)   {\cal A} L_k(x)  {dx \over x}.\eqno(3.31)$$ 
However,  taking into account the second order linear differential equation for Laguerre polynomials [3, formula (4.18.8)], we see that the action of the operator (1.5) on  the Laguerre polynomial $L_k(x)$ gives  the relation 

$${\cal A} L_k(x) = x^2 L_k(x) + k x L_{k-1}(x).$$ 
Hence due to the 3-term recurrence relation for Laguerre polynomials we find from (3.19), (3.31)

$$l_{m,k} =  -  \int_0^\infty e^{-x}  p_{m-1} (x) x  L_k(x)  dx -  k \int_0^\infty e^{-x}  p_{m-1} (x)   L_{k-1} (x)  dx $$

$$= k(k-1) l_{m-1,k}- k(2k-1) l_{m-1,k-1} + k(k-1)  l_{m-1,k-2}$$

$$ + (k+1) l_{m-1,k+1} - (2k+1) l_{m-1,k} + k l_{m-1,k-1},$$
i.e. after simplification

$$l_{m,k} = (k+1) l_{m-1,k+1} + (k^2- 3k-1) l_{m-1,k} + 2k( 1-k) l_{m-1,k-1} $$

$$+ k(k-1)  l_{m-1,k-2},\quad  m,k \in \mathbb{N}.\eqno(3.32)$$
In the meantime, returning to (3.15) and invoking (3.25), we get

$$H_{n+1} c_{n+1,k} + E_n   c_{n,k}  + H_n  c_{n-1,k} = \int_0^\infty  e^{-x} L_k(x) {\cal A}  F_n(x) dx$$

$$=  \int_0^\infty  {\cal A}  \left[ e^{-x} x L_k(x) \right] F_n(x) {dx\over x} .$$
Recalling (1.5) and the  differential equation for Laguerre polynomials, we have, in turn,

$$ {\cal A}  \left[ e^{-x} x L_k(x) \right] =  e^{-x} x^3 L_k(x) -   x{d\over dx} x {d\over dx} \left[ e^{-x} x L_k(x) \right] $$

$$= e^{-x} x^3 L_k(x) -   x{d\over dx} \left[ - e^{-x} x^2 L_k(x) + e^{-x} x L_k(x) + e^{-x} x^2 L^\prime_k(x) \right] $$ 

$$=  -   e^{-x} x  \left[  L_k(x) +  3  x \left( L^\prime_k(x) - L_k(x) \right) + x^2 \left( L^{\prime\prime}_k(x) - 2L^\prime_k(x)\right)\right] $$ 

$$=    e^{-x} x \bigg[  \left( (3+2k) x -2k-1\right) L_k(x)  +   (2-x) k\ L_{k-1} (x) \bigg].$$ 
Consequently, with the 3-term recurrence relation for Laguerre polynomials we derive 

$$H_{n+1} c_{n+1,k} + E_n   c_{n,k}  + H_n  c_{n-1,k} = 2k c_{n,k-1}  - (2k+1)  c_{n,k} +  (3+2k) \int_0^\infty e^{-x} x L_k(x) F_n(x)  dx$$

$$  - k \int_0^\infty e^{-x} x  L_{k-1} (x) F_n(x)  dx  $$

$$= 2k c_{n,k-1}  - (2k+1)  c_{n,k} +  (3+2k) \int_0^\infty e^{-x} \left[ (2k+1) L_k(x) - (k+1)L_{k+1} (x) - k L_{k-1} (x) \right] F_n(x)  dx$$

$$  + k \int_0^\infty e^{-x} \left[ k L_{k} (x) - (2k-1) L_{k-1} (x) + (k-1) L_{k-2} (x) \right] F_n(x)  dx  $$

$$= - (k+1)(3+2k) c_{n,k+1} + ( 5k^2 + 6k+2) c_{n,k} - 4k^2 c_{n,k-1} + k(k-1) c_{n,k-2}.$$  
Thus we establish, finally, the following identity

$$ H_{n+1} c_{n+1,k} + E_n   c_{n,k}  + H_n  c_{n-1,k} $$

$$= - (k+1)(3+2k) c_{n,k+1} + ( 5k^2 + 6k+2) c_{n,k} - 4k^2 c_{n,k-1} + k(k-1) c_{n,k-2},\eqno(3.33)$$ 
where (see above) $n, k  \in \mathbb{N}_0,\   c_{n, -2} = c_{n, -1} = c_{n, 0} =\dots= c_{n, n-2} =0.$ Letting $k=n-2$,  we get immediately (3.30).  One can express the coefficients $c_{n,k},\ k\ge n-1$ in terms of determinants.  Indeed,  recalling (3.22) and using (3.1), we find for $m > n$

$$ \sum_{k=n}^{m}  (-1)^k f_{m,k} \sum_{q=n}^{k}  c_{n,q-1} l_{k,q-1} = \sum_{q=n}^{m}  c_{n,q-1} \sum_{k=q}^{m}  (-1)^k f_{m,k}\ l_{k,q-1} = 0.$$
Denoting by 

$$s_{m,q} = \sum_{k=q}^{m}  (-1)^k f_{m,k}\ l_{k,q-1}$$
we write the non-homogeneous linear system of $m$ equations with $m$ unknowns 

$$c_{n,n}, c_{n,n+1},\dots, c_{n,n+m-1},\ m \in \mathbb{N}$$
 with nonzero determinant of a lower triangular matrix. In fact, having in mind (3.23),  this gives 
 
 $$  \sum_{q=n+1}^{n+\nu}  c_{n,q-1} s_{n+\nu,q} =  {(-2)^{n-1} \pi  s_{n+\nu,n} \over f_n\  (2n-1)! },\ \nu = 1,2,\dots, m,$$
where the corresponding determinant has the value $\prod_{j=1}^m s_{n+j, n+j}$ and the unique solution takes the form, accordingly,  due to the Cramer rule.

Further,  regarding (3.13) as the inverse Kontorovich-Lebedev transform [6, Ch. 2],  this means, reciprocally,  

$$ W_n(\tau^2, [a]_{r+1}) = {2\over \pi^2} \bigg|\frac{\Gamma(2i\tau)} { \prod_{k=1}^{r+1}  \Gamma\left( a_k + i\tau\right)}\bigg|^{2}  \int_0^\infty K_{i\tau} (x)  F_n(x) {dx\over x},\eqno(3.34)$$
where the integral (3.34) converges at least in the mean square sense with respect to the norm in $L_2 (\mathbb{R}_+; \tau\sinh(\pi\tau) d\tau)$.   Moreover, the Parseval equality holds

$$ \int_0^\infty  F^2 _n(x) {dx\over x}  = {\pi^2\over 2} \int_0^\infty \tau\sinh(\pi\tau) \bigg|\frac { \prod_{k=1}^{r+1}  \Gamma\left( a_k + i\tau\right)} {\Gamma(2i\tau)} \bigg|^{4} W^2_n(\tau^2, [a]_{r+1})d\tau.\eqno(3.35)$$
Functions $F_n(x)$ can be expressed in the operator form as polynomials $W_n$. Indeed, taking into account (1.5), (1.6), (3.10), (3.13) and termwise differentiation  under the integral sign, we have

$$F_n(x)= \sum_{k=0}^n f_{n,k}  \int_0^\infty \tau \sinh(\pi\tau)  \bigg|\frac{ \prod_{k=1}^{r+1}  \Gamma\left( a_k + i\tau\right)}{\Gamma(2i\tau)} \bigg|^2    K_{i\tau} (x) \tau^{2k} d\tau$$

$$= \sum_{k=0}^n f_{n,k} {\cal A}^k   \int_0^\infty \tau \sinh(\pi\tau)  \bigg|\frac{ \prod_{k=1}^{r+1}  \Gamma\left( a_k + i\tau\right)}{\Gamma(2i\tau)} \bigg|^2    K_{i\tau} (x)  d\tau$$

$$= W_n({\cal A}, [a]_{r+1}) \left[ \Phi(x)\right],\eqno(3.36)$$
where $\Phi(x)$ is the following Kontorovich-Lebedev integral

$$\Phi(x)= \int_0^\infty \tau \sinh(\pi\tau)  \bigg|\frac{ \prod_{k=1}^{r+1}  \Gamma\left( a_k + i\tau\right)}{\Gamma(2i\tau)} \bigg|^2    K_{i\tau} (x)  d\tau,\quad x >0.\eqno(3.37)$$
Reciprocally, we get the equality (compare with (3.34))
$$ {\pi^2\over 2} \bigg|\frac{ \prod_{k=1}^{r+1}  \Gamma\left( a_k + i\tau\right)}{\Gamma(2i\tau)} \bigg|^2 =  \int_0^\infty  K_{i\tau} (x) \Phi(x) {dx\over x}.$$
Besides, recalling (1.19), (3.36), we rewrite (3.12) in the form
$$  \int_0^\infty e^{-x}  p_{m+1} (x)  F_n(x) {dx \over x} =  (-1)^{m+1} \int_0^\infty   {\cal A} W_n({\cal A}, [a]_{r+1}) \left[ e^{-x} \right] $$

$$\times  {\cal A}^m \left[ \Phi(x)\right] {dx \over x}  = 0,\quad   m=0,1,2,\dots, n-2.\eqno(3.38)$$
Therefore the Parseval equality for the Kontorovich-Lebedev transform yields

$$ \int_0^\infty \tau \sinh(\pi\tau)  \bigg|\frac{ \prod_{k=1}^{r+1}  \Gamma\left( a_k + i\tau\right)}{\Gamma(2i\tau)} \bigg|^2 \tau^{2m} $$

$$\times \int_0^\infty  {\cal A}  W_n({\cal A}, [a]_{r+1}) \left[ e^{-x} \right]   K_{i\tau} (x)  {dx d\tau\over x}  = 0,\  m=0,1,2,\dots, n-2.\eqno(3.39)$$
Meanwhile, via (1.19) and (1.26) we find

$$ \tau \sinh(\pi\tau) \int_0^\infty  {\cal A}  W_n({\cal A}, [a]_{r+1}) \left[ e^{-x} \right]   K_{i\tau} (x)  {dx\over x} $$

$$= \tau \sinh(\pi\tau) \sum_{k=0}^n f_{n,k} (-1)^{k+1} \int_0^\infty   e^{-x} p_{k+1} (x)  K_{i\tau} (x)  {dx\over x}$$

$$ = \tau \sinh(\pi\tau) \sum_{k=0}^n f_{n,k} (-1)^{k} \sum_{\nu =0}^k \binom{2k+1}{2\nu} \int_0^\infty   e^{-x} p_{\nu} (x)  K_{i\tau} (x) dx.\eqno(3.40)$$
Moreover,  Entry 3.14.3.1 in [1], (1.5), (1.6), (3.10)  and the property (1.35) of the Macdonald function show that the integral 

$${1\over \pi} \ \tau \sinh(\pi\tau)  \int_0^\infty  {\cal A}  W_n({\cal A}, [a]_{r+1}) \left[ e^{-x} \right]   K_{i\tau} (x)  {dx \over x} $$
is equal to $\tau^2 W_n(\tau^2, [a]_{r+1})$.  Precisely, owing to the integration by parts,  we derive 

$${1\over \pi} \ \tau \sinh(\pi\tau)  \int_0^\infty  {\cal A}  W_n({\cal A}, [a]_{r+1}) \left[ e^{-x} \right]   K_{i\tau} (x)  {dx d\tau\over x} $$

$$= {1\over \pi} \ \tau^2 \sinh(\pi\tau) \sum_{k=0}^n f_{n,k}   \tau^{2k} \int_0^\infty  e^{-x} {\rm Im} \left[ K_{1+i\tau}(x) \right]   dx $$

$$= \sum_{k=0}^n f_{n,k}   \tau^{2(k+1)} = \tau^2 W_n(\tau^2, [a]_{r+1}).$$
Now, returning to (3.13), we multiply both sides by $e^{-x} L_k(x)$ and integrate over  $\mathbb{R}_+$. Then, interchanging the order of integration on the right-hand side, owing to the dominated convergence, we find the values of the coefficients $c_{n,k}$ in the form

$$c_{n,k} =  \int_0^\infty \tau \sinh(\pi\tau)  \bigg|\frac{ \prod_{m=1}^{r+1}  \Gamma\left( a_m + i\tau\right)}{\Gamma(2i\tau)} \bigg|^2 W_n(\tau^2, [a]_{r+1}) $$

$$\times   \int_0^\infty  e^{-x} L_k(x)  K_{i\tau} (x)  dx d\tau.\eqno(3.41)$$
The inner integral by $x$ is calculated in [1, Entry 3.23.11.3] in terms of the generalized hypergeometric functions ${}_4F_3$ at the unity, and we have the formula 

$$ \int_0^\infty  e^{-x} L_k(x)  K_{i\tau} (x)  dx  = {\pi\over 2\cosh(\pi\tau/2)} \  {}_4F_3\left({1+k\over 2},\  1+{k\over 2},\ {1+i\tau\over 2}, \ {1-i\tau\over 2}; \ {1\over 2}, {1\over 2}, 1;  1\right) $$

$$+ {(k+1)\pi \tau\over 2\sinh(\pi\tau/2)} \  {}_4F_3\left({3+k\over 2},\  1+{k\over 2},\  1+ {i\tau\over 2}, \  1- {i\tau\over 2}; \ {3\over 2}, {3\over 2}, 1;  1\right).$$
Thus we derive the values $c_{n,k} $ in terms of the integrals, involving polynomials $W_n(\tau^2, [a]_{r+1})$

$$ c_{n,k} =  \pi \int_0^\infty \tau \sinh\left({\pi\tau\over 2}\right)  \bigg|\frac{ \prod_{m=1}^{r+1}  \Gamma\left( a_m + i\tau\right)}{\Gamma(2i\tau)}  \bigg|^2  W_n(\tau^2, [a]_{r+1}) $$

$$\times {}_4F_3\left({1+k\over 2},\  1+{k\over 2},\ {1+i\tau\over 2}, \ {1-i\tau\over 2}; \ {1\over 2}, {1\over 2}, 1;  1\right)  d\tau $$

$$+ \pi (k+1) \int_0^\infty   \cosh\left({\pi\tau\over 2}\right)  \bigg|\frac{ \prod_{m=1}^{r+1}  \Gamma\left( a_m + i\tau\right)}{\Gamma(2i\tau)}  \bigg|^2   \tau^2 W_n(\tau^2, [a]_{r+1}) $$

$$\times    {}_4F_3\left({3+k\over 2},\  1+{k\over 2},\  1+ {i\tau\over 2}, \  1- {i\tau\over 2}; \ {3\over 2}, {3\over 2}, 1;  1\right) d\tau,\ (n,k) \in \mathbb{N}_0\times\mathbb{N}_0.$$
Now, returning to (3.34), we write

$$ {\pi^2\over 2}  \bigg|\frac { \prod_{k=1}^{r+1}  \Gamma\left( a_k + i\tau\right)} {\Gamma(2i\tau)}\bigg|^{2} W_n(\tau^2, [a]_{r+1}) =   \int_0^\infty K_{i\tau} (x)  F_n(x) {dx\over x}.$$
Hence we multiply both sides of the latter equality by  $W_m(\tau^2, [a]_{r+1})$ and integrate over $\mathbb{R}_+$. Then, taking in mind (3.1) and (3.34), we interchange the order of integration which can be justified by the Cauchy-Schwarz inequality and dominated convergence theorem to obtain

$$ \int_0^\infty   \int_0^\infty  W_m(\tau^2, [a]_{r+1}) K_{i\tau} (x)  F_n(x) { d\tau  dx\over x} = {\pi^2\over 2}\delta_{n,m}.$$
Meanwhile,  equalities (1.1), (1.6), (3.36) and  integration by parts in the inner integral by $ x$ yield  the following composition orthogonality 

$$ \int_0^\infty   \int_0^\infty  W_m(\tau^2, [a]_{r+1}) K_{i\tau} (x)  F_n(x) { d\tau  dx\over x} $$

$$=   \int_0^\infty   \int_0^\infty  W_m({\cal A}, [a]_{r+1}) W_n({\cal A}, [a]_{r+1}) \left[ K_{i\tau} (x) \right]   \Phi(x)  { d\tau  dx\over x} $$

$$=      \int_0^\infty  W_m({\cal A}, [a]_{r+1}) W_n({\cal A}, [a]_{r+1}) \left[ e^{-x} \right]   \Phi(x)  {dx\over x}  =  \pi \delta_{n,m}.$$
On the other hand,  employing in (3.41) the Rodrigues formula for Laguerre polynomials, integration by parts and the identity for derivatives of the Macdonald function [4, Vol. II]

$${d^k\over dx^k} K_{i\tau}(x)= {(-1)^k\over 2^k} \sum_{\nu=0}^k \binom{k}{\nu} K_{2\nu-k+i\tau} (x),$$ 
we have 

$$c_{n,k} =  {1\over 2^k  k!}  \sum_{\nu=0}^k \binom{k}{\nu} \int_0^\infty \tau \sinh(\pi\tau)  \bigg|\frac{ \prod_{m=1}^{r+1}  \Gamma\left( a_m + i\tau\right)}{\Gamma(2i\tau)} \bigg|^2 W_n(\tau^2, [a]_{r+1}) $$

$$\times   \int_0^\infty  e^{-x} x^k  K_{2\nu-k+i\tau} (x)  dx d\tau.$$
The latter inner integral by $x$ is calculated via Entry  3.14.3.1 in [1] to obtain the expression after simplification

$$c_{n,k} =  {\pi \over (2k+1)!}  \sum_{\nu=0}^k \binom{k}{\nu} \int_0^\infty  \bigg|\frac{ \prod_{m=1}^{r+1}  \Gamma\left( a_m + i\tau\right)}{\Gamma(2i\tau)} \bigg|^2   W_n(\tau^2, [a]_{r+1})$$

$$\times  (-i\tau)_{2(k-\nu)+1} (i\tau)_{2\nu+1} d\tau.\eqno(3.42)$$
Decomposing Pochhammer symbols in terms of the Stirling numbers of the first kind $s(n,k)$ [11], it gives 

$$(-i\tau)_{2(k-\nu)+1} = \sum_{\lambda=0}^{2(k-\nu)+1} (-1)^{2(k-\nu)+1-\lambda} s(2(k-\nu)+1,\lambda) (-i\tau)^\lambda,$$

$$(i\tau)_{2\nu+1} = \sum_{\omega=0}^{2\nu+1} (-1)^{2\nu+1-\omega} s(2\nu+1,\omega) (i\tau)^\omega.$$
Hence the finite real sum  in (3.42) can be treated as follows 

$$S_{k+1}(\tau^2)= {\pi \over (2k+1)!}  \sum_{\nu=0}^k \binom{k}{\nu}  (-i\tau)_{2(k-\nu)+1} (i\tau)_{2\nu+1}$$

$$ = {\pi \over (2k+1)!}  \sum_{\nu=0}^k \binom{k}{\nu}  {\rm Re} \bigg[  (-i\tau)_{2(k-\nu)+1} (i\tau)_{2\nu+1}\bigg]$$

$$= {\pi \over (2k+1)!}  \sum_{\nu=0}^k \sum_{\lambda=0}^{k-\nu} \sum_{\omega=0}^{\nu}  (-1)^{\lambda+\omega}  \binom{k}{\nu}  \tau^{2(\lambda+\omega)} \bigg[ s(2(k-\nu)+1,2\lambda)\  s(2\nu+1, 2\omega)\bigg.$$

$$\bigg.  + \tau^2  s(2(k-\nu)+1,2\lambda+1)\  s(2\nu+1, 2\omega+1)\bigg],\ k \in \mathbb{N}_0,$$
and we see that $S_{k+1}(\tau^2)$  represents a polynomial of degree $k+1$ in $\tau^2$.  Consequently,  the orthogonality (3.1) ensures the polynomial expansion

$$S_{k+1}(\tau^2) = \sum_{j=0}^{k+1} c_{j,k} W_j(\tau^2, [a]_{r+1})$$
and the Parseval equality

$$\sum_{j=0}^{k+1} c^2_{j,k} =  \int_0^\infty  \bigg|\frac{ \prod_{m=1}^{r+1}  \Gamma\left( a_m + i\tau\right)}{\Gamma(2i\tau)} \bigg|^2 S^2_{k+1}(\tau^2) d\tau.$$

\bigskip
\centerline{{\bf Acknowledgments}}
\bigskip

The work was partially supported by CMUP, which is financed by national funds through FCT (Portugal) under the project with reference UID/00144/2025.

\bigskip
\centerline{{\bf References}}
\bigskip
\baselineskip=12pt
\medskip
\begin{enumerate}

 \item[{\bf 1.}\ ]
Yu. A. Brychkov,  O.I. Marichev,  N.V.  Savischenko, {\it Handbook of Mellin transforms. Advances in Applied Mathematics},  CRC Press, Boca Raton, FL, 2019. 

 \item[{\bf 2.}\ ]
  R. Koekoek, P.A. Lesky, R.F. Swarttouw, {\it Hypergeometric Orthogonal Polynomials and their q-Analogues}, in: Springer Monographs in Mathematics,
Springer Verlag, New York, 2010.

\item[{\bf 3.}\ ]
 N.N. Lebedev, {\it Special Functions and Their Applications}, Dover,  New York, 1972.

\item[{\bf 4.}\ ] A.P. Prudnikov, Yu.A. Brychkov and O.I. Marichev, {\it Integrals and Series}. Vol. I: {\it Elementary Functions}, Vol. II: {\it Special Functions}, Gordon and Breach, New York and London, 1986, Vol. III: {\it More Special Functions}, Gordon and Breach, New York and London, 1990.

\item[{\bf 5.}\ ]S. Yakubovich and Yu.  Luchko, {\it The Hypergeometric Approach to Integral Transforms and Convolutions},
(Kluwers Ser. Math. and Appl.: Vol. 287), Dordrecht, Boston, London, 1994.

\item[{\bf 6.}\ ] S. Yakubovich, {\it Index Transforms}, World Scientific Publishing Company, Singapore, New Jersey, London and Hong Kong, 1996.

\item[{\bf 7.}\ ] S. Yakubovich, Certain identities, connection and explicit formulas for the Bernoulli and Euler numbers and the Riemann zeta-values, {\it Analysis} {\bf 35} (2015), 1,  59- 71.

\item[{\bf 8.}\ ] S. Yakubovich,  On the Yor integral and a system of polynomials related to the Kontorovich-Lebedev transform, {\it Integral Transforms and Special Functions} {\bf 24} (2013), 8,  672-683.

\item[{\bf 9.}\ ] S. Yakubovich, On the orthogonality and convolution orthogonality via the Kontorovich-Lebedev transform, {\it Journal of Computational and Applied Mathematics}, {\bf 384} (2021),  paper No. 113178.

\item[{\bf 10.}\ ] S. Yakubovich, A method of composition orthogonality and new sequences of orthogonal polynomials and functions for non-classical weights, {\it J. Math. Anal. Appl. } {\bf 499}  (2021), N 2, paper No. 125032.

\item[{\bf 11.}\ ] E.W. Weisstein, Pochhammer Symbol. From MathWorld- A Wolfram Resource. https://mathworld.wolfram.com/PochhammerSymbol.html

\end{enumerate}

\vspace{5mm}

\noindent Semyon  Yakubovich\\
Department of  Mathematics,\\
Faculty of Sciences,\\
University of Porto,\\
Campo Alegre st., 687\\
4169-007 Porto\\
Portugal\\
E-Mail: syakubov@fc.up.pt\\

\end{document}